\documentclass[11pt,a4paper,reqno,centertags]{article}
\usepackage{xypic,a4wide,amsmath,amsthm,amssymb}

\newtheorem{thm}{Theorem}[section]
\newtheorem{corol}[thm]{Corollary}
\newtheorem{lemma}[thm]{Lemma}
\newtheorem{prop}[thm]{Proposition}
\newtheorem{defin}[thm]{Definition}

\theoremstyle{remark}
\newtheorem{rem}[thm]{Remark}
\newtheorem{ex}[thm]{Example}

\newcommand{\qee}{\mbox{\hspace{0.2mm}}\hfill$\triangle$}
\def\prf{\noindent{\textsc{Proof}}\rm\ }
\def\endprf{\ \hfill $\Box$\medskip}

\def\cal{\mathcal}

\def\bC{{\Bbb C}}

\def\bP{{\Bbb P}}

\def\bZ{{\Bbb Z}}
\def\bN{{\Bbb N}}

\def\unmez{\frac{1}{2}}

{\nonumber}

\begin{document}
\begin{center}
 {\LARGE\bf Block Toeplitz determinants, constrained $\mathrm{KP}$ and Gelfand-Dickey hierarchies.}\\[15pt]
 {\sc Mattia Cafasso}\\
 {(SISSA - Trieste)}
\end{center}

\vspace{1cm}

\begin{abstract}
	We propose a method for computing any Gelfand-Dickey $\tau$ function living in Segal-Wilson Grassmannian as the asymptotics of block Toeplitz determinant associated to a certain class of symbols $\cal W(t;z)$.
	Also truncated block Toeplitz determinants associated to the same symbols are shown to be $\tau$ function for rational reductions of $\mathrm{KP}$.\\
	Connection with Riemann-Hilbert problems is investigated both from the point of view of integrable systems and block Toeplitz operator theory.\\
	Examples of applications to algebro-geometric solutions are given.
\end{abstract}

\section*{Introduction}
  This paper deals with the applications of block Toeplitz determinants and their asymptotics to the study of integrable hierarchies.
	Asymptotics of block Toeplitz determinants and their applications to physics is a developing field of research;
	in recent years it has been shown how to compute some physically relevant quantities (e.g. correlation functions) studying asymptotics of some block Toeplitz determinants (see \cite{IJK},\cite{IMM},\cite{BE}).
	In particular in \cite{IJK} and \cite{IMM} authors, for the first time, showed effective computations for the case of block Toeplitz determinants with symbols that do not have half truncated Fourier series. This is of particular interest for us as, with our approach, we will be able to do the same for certain block Toeplitz determinants associated to algebro-geometric solutions of Gelfand Dickey hierarchies. 
	Let us mention some theoretical results about (block) Toeplitz determinants we will use in this paper. 
	Given a function $\gamma(z)$ on the circle we denote $T_N(\gamma)$ the Toeplitz matrix with symbol $\gamma$ given by
$$
T_N(\gamma):=
\begin{pmatrix}
\gamma^{(0)} & \ldots & \ldots & \gamma^{(-N)}\\
&&&\\
\gamma^{(1)} & \ldots & \ldots & \gamma^{(-N+1)}\\
&&&\\
   \ldots       & \ldots & \ldots & \ldots\\
&&&\\
\gamma^{(N)} & \ldots & \ldots & \gamma^{(0)}
\end{pmatrix}$$
where $\gamma^{(k)}$ are the Fourier coefficients $\gamma(z)=\sum_{k} \gamma^{(k)}z^k$.\\
We use the term \emph{block Toeplitz} for the case of matrix-valued symbol $\gamma(z)$. In that case the entries $\gamma^{(j-i)}$ of the above matrix are $n\times n$ matrices themselves. We denote 
$$D_N(\gamma):=\det T_N(\gamma)$$ 
and we use the notation $T(\gamma)$ for the $\bN\times\bN$ matrix obtained letting $N$ go to infinity. The main goal of the theory of Toeplitz determinants is to compute $D_N(\gamma)$ as $N$ goes to infinity and find expressions for $D_N(\gamma)$ as well as for its limit in terms of Fredholm determinants.\\
First result is due to Szeg\"o that in 1952 gave a formula for asymptotics of $D_N(\gamma)$ in the scalar case \cite{Sz}. This result has been generalized by H. Widom in the 70's (\cite{W0},\cite{W1} and \cite{W2}) for the matrix case; namely he proved that under suitable analytical assumptions it exists the limit 
		$$D_\infty(\gamma):=\lim_{N\rightarrow\infty} \frac{D_N(\gamma)}{G(\gamma)^N}=\det(T(\gamma)T(\gamma^{-1}))$$
where $G(\gamma)$ is a normalizing constant and the operator $T(\gamma)T(\gamma^{-1})$ is such that its determinant is well defined as a Fredholm determinant (see section 2 for the precise statement).
Once the asymptotics had been computed the next quite natural question was to find an expressions relating directly $D_N(\gamma)$, and not just its asymptotics, to certain Fredholm determinants. The problem was solved many years later by Borodin and Okounkov in  \cite{BO} for the scalar case and generalized, in the same year, for matrix case by E. Basor and H.Widom in $\cite{BW}$.
For matrix valued case Borodin-Okounkov formula reads 
$$D_N(\gamma)=D_\infty(\gamma)\det(I-K_{\gamma,N})$$
(here we assume $G(\gamma)=1$).
The operator $(I-K_{\gamma,N})$ can be written explicitely in coordinates knowing certain Riemann-Hilbert factorizations of $\gamma$. Its Fredholm determinant is well defined (see section 2 for details).
Now many proofs of Borodin-Okounkov formula are known (for instance \cite{Bot} contains another proof of the same formula, see also the earlier paper \cite{GC}).\\

In this paper we apply block-Toeplitz determinants to the computation of $\tau$ function of an (almost) arbitrary solution of Gelfand-Dickey hierarchy
	$$\frac{\partial L}{\partial t_j}=[(L^{\frac{j}{n}})_+,L].$$
	($L$ differential operator of order $n$, $j\neq nk$).\\
	 More precisely to a given point 
	 $$W=\cal W(z)H^{(n)}_+$$
	 in the big cell of Segal-Wilson vector-valued Grassmannian we associate a $n\times n$ matrix-valued symbol $\cal W(t;z)$ obtained deforming $\cal W(z)$ (see formula (\ref{GDSymbols})).
	 In this way we define a sequence of $N$-truncated block Toeplitz determinants $\lbrace\tau_{W,N}(t)\rbrace_{N>0}$ which are shown to be solutions of certain rational reductions of $\mathrm{KP}$; this is our\\
	 
	 {\bf First result:} \emph{Every symbol $\cal W(t;z)$ defines through its truncated determinants a sequence $\lbrace\tau_{W,N}(t)\rbrace_{N>0}$ of solutions for $\mathrm{KP}$ such that}
	 $$\tau_{W,N}(t)\in\mathrm{cKP}_{1,nN}\cap\mathrm{cKP}_{n,n}\quad\forall N>0.$$
	 Here we used the notation from \cite{Ar}; given a $\tau$ function for $\mathrm{KP}$ with corresponding Lax pseudodifferential operator $\cal L$ we say that $\tau\in\mathrm{cKP}_{m,n}$ iff $\cal L^m$ can be written as the ratio of two differential operators of order $m+n$ and $n$ respectively.
	 This sequence admits a \emph{stable limit} which is shown to be equal to the Gelfand-Dickey $\tau$ function $\tau_W(t)$ associated to $W$; this quantity can be computed using Szeg\"o-Widom's theorem. This will give us the remarkable identity
	 \begin{equation}\label{introW}
	 		\tau_{W}(\tilde t)=\det\Big[\cal P_{\cal W(\tilde t;z)}\Big]
	 \end{equation}
	 where $\cal P_{W(\tilde t;z)}$ is the Fredholm operator appearing in Szeg\"o-Widom's theorem (here we put $\tilde t$ instead of $t$ to remember that, when working with $W\in \mathrm{Gr^{(n)}}$, times $t_{nj}$ multiple of $n$ must be set to $0$). Next step is the study of Riemann-Hilbert (also called Wiener-Hopf) factorization of symbol $\cal W(t;z)$ given by
	 \begin{equation}\label{introRH}
	 		\cal W(t;z)=T_-(t;z)T_+(t;z)
	 \end{equation}
	 with $T_-$ and $T_+$ analytical in $z$ outside and inside $S^1$ respectively and normalized as 
	 $$T_-(\infty)=I.$$
	 Here we assume that the symbol can be extended to an analytic function in a neighborhood of $S^1$.
	 Using Plemelj's work \cite{P} we show that $T_-(\tilde t;z)$ must satify the integral equation
	 \begin{equation}\label{introintegralequation}
	 		\cal P_{\cal W(\tilde t;z)}^{\mathrm T}T_-(\tilde t;z)=I
	 \end{equation}
	 and  we write a solution of (\ref{introRH}) in terms of wave function $\psi_W(\tilde t;z)$ corresponding to $W$.\\
	 In this way we arrive to our second result:\\
	 
	 {\bf Second result:}\emph{ Take $W\in\mathrm{Gr}^{(n)}$ in the big cell and its corresponding $\tau$ function $\tau_{W}(\tilde t)$.\\ 
	 $\tau_W(\tilde t)$ is equal to the Fredholm determinant of the homogeneous integral equation associated to (\ref{introintegralequation}) which is related to Riemann-Hilbert problem (\ref{introRH}).
	 The solution of this Riemann-Hilbert problem is unique for every value of parameters $\tilde t$ that makes $\tau_W(\tilde t)\neq 0$ and can be computed by means of related wave function $\psi_W(\tilde t;z)$.}\\

At the end of the paper we consider a particular class of symbols $\cal W(t;z)$ corresponding to algebro-geometric solutions of Gelfand-Dickey hierarchies. We formulate an alternative Riemann-Hilbert problem equivalent to (\ref{introRH}) and explain how to solve it using $\theta$-functions. In this way we give concrete formulas for a wide class of symbols that do not have half truncated Fourier series. We think this is quite remarkable since concrete results for non half truncated symbols were available, till now, just for the concrete cases presented in \cite{IJK} and \cite{IMM}.\\
The paper is organized as follows:
\begin{itemize}
	\item First section states some results about Segal-Wilson Grassmannian and related loop groups we will need in the sequel; proofs can be found in \cite{SW} and \cite{PS}.
	\item Second section states Szeg\"o-Widom's theorem and related results obtained by Widom in \cite{W0},\cite{W1} and \cite{W2} and the Borodin-Okounkov formula for block Toeplitz determinant \cite{BW}.
	\item In the third section we introduce and study the sequence of truncated determinants $\lbrace\tau_{W,N}(t)\rbrace_{N>0}$ and its stable limit $\tau_W(t)$. We want to remark that the property of stability was stated for the first time in \cite{IZ} (see also \cite{DF}) and our sequence is actually a subsequence of the stabilizing chain studied in \cite{D1}; nevertheless, to our best knowledge, this is the first time that block Toeplitz determinants enter the game and also the observation that $\tau_{W,N}\in\mathrm{cKP}_{n,n}$ seems to be something new.
	\item Fourth section is devoted to establishing the connection between integral equations formulated by Plemelj in \cite{P} and Fredholm operator appearing in Szeg\"o-Widom's theorem.
	\item In the fifth section we show how to write Riemann-Hilbert factorization of $\cal W(\tilde t;z)$ in terms of wave function $\psi_W(\tilde t;z)$. Of course relation between Gelfand-Dickey hierarchy and factorization problem is something known; our exposition here is closely related to \cite{SS}. Moreover, knowing Riemann-Hilbert factorization of $\cal W(\tilde t;z)$, we can apply Borodin-Okounkov formula to give an expression of any $\tau_{W,N}(\tilde t)$ as Fredholm determinant and a recursion relation to go from $\tau_{W,N}(\tilde t)$ to $\tau_{W,N+1}(\tilde t)$.
	\item Last section gives explicit formulas for symbols and $\tau$ functions associated to algebro-geometric rank one solutions of Gelfand-Dickey hierarchies. Also we formulate an alternative Riemann-Hilbert problem equivalent to (\ref{introRH}) in analogy with what has been done in \cite{IJK} and \cite{IMM}. We explain how to solve it using $\theta$-functions.
	
\end{itemize}

\subsection*{Acknowledgements} 

I am grateful to A.Its for a fruitful discussion in which he pointed out that our Riemann-Hilbert problem can be reduced to a problem with constant jump as in \cite{IJK},\cite{IMM} and gave me some interesting hints about possible developments of this work.\\
J. van de Leur suggested me to verify if vector-costrained reductions of $\mathrm{KP}$ had some relations with block Toeplitz determinant, I am grateful to him for this really useful suggestion.\\
Moreover I wish to thank my advisor B.Dubrovin for his constant support and suggestions he gave me during many hours spent together discussing the preparation of this paper. 

\vspace{0.5cm}
 
This work is partially supported by the European Science Foundation Programme ``Methods of Integrable Systems, Geometry, Applied Mathematics" (MISGAM), the Marie Curie RTN ``European Network in Geometry, Mathematical Physics and Applications"  (ENIGMA),  and by the Italian Ministry of Universities and Researches (MUR) research grant PRIN 2006 ``Geometric methods in the theory of nonlinear waves and their applications".

\section{Segal-Wilson Grassmannian and related loop groups}
Here we recall some definitions and results from \cite{SW} and \cite{PS} that will be useful in the sequel.
\begin{defin}
	Let $H^{(n)}:=L^2(S^1,\bC^n)$ be the space of complex vector-valued square-integrable functions.
	We choose a orthonormal basis given by 
	$$\lbrace e_{\alpha,k}:=(0,\ldots,z^k,\ldots,0)^{\mathrm{T}}:\alpha=1\ldots n, k\in\bZ\rbrace $$
	and the polarization
  $$H^{(n)}=H^{(n)}_+\oplus H^{(n)}_-$$
  where $H^{(n)}_+$ and $H^{(n)}_-$ are the closed subspaces spanned by elements $\{e_{\alpha,k}\}$ with $k\geq 0$ and $k<0$ 		  respectively.
\end{defin}
In the sequel in order to avoid cumbersome notations we will write $H$ instead of $H^{(1)}$.
\begin{defin}[\cite{PS}]
   The Grassmannian $\mathrm{Gr}(H^{(n)})$ modeled on $H^{(n)}$ consists of the subset of closed subspaces $W\subseteq H^{(n)}$ such that:
	\begin{itemize}
		\item the orthogonal projection $\mathrm{pr}_+:W\rightarrow H^{(n)}_+$ is a Fredholm operator.
		\item the orthogonal projection $\mathrm{pr}_-:W\rightarrow H^{(n)}_-$ is a Hilbert-Schmidt operator.
	\end{itemize}
	Moreover we will denote $\mathrm{Gr}^{(n)}$ the subset of $\mathrm{Gr}(H^{(n)})$ given by subspaces $W$ such that $zW\subseteq W$.
\end{defin}
It's well known \cite{SW} that through Segal-Wilson theory we can associate a solution of $n^{th}$ Gelfand-Dickey hierarchy to every element of $\mathrm{Gr}^{(n)}$; this is the reason why we are interested in them.
\begin{lemma}[\cite{SW}]
   The map 
   \begin{align*}
      \Xi:H^{(n)} & \longrightarrow H\\
        (f_0(z),\ldots,f_{n-1}(z))^{\mathrm{T}}  & \longmapsto\ \tilde{f}(z):=f_0(z^n)+\ldots+z^{n-1}f_{n-1}(z^n)
		\end{align*}
		is an isometry. Its inverse is given by 
		$$f_k(z)=\frac{1}{n}\sum_{\zeta^n=z}\zeta^{-k}\tilde{f}(\zeta)$$
		where the sum runs over the $n^{th}$ roots of $z$.
\end{lemma}
\begin{prop}
     Under the isometry $\Xi$ we can identify $\mathrm{Gr}^{(n)}$ with the subset 
     $$\lbrace W\in \mathrm{Gr}(H):z^nW\subseteq W\rbrace$$ 
\end{prop}

It is obvious that loop groups act on Hilbert spaces defined above by multiplications.
We want to define a certain loop group $L_{\unmez}\mathrm{Gl}(n,\bC)$ with good analytical properties acting transitively on $\mathrm{Gr}^{(n)}$; in such a way we can obtain any $W\in \mathrm{Gr}^{(n)}$ just acting on the reference point $H^{(n)}_+$ with this group.
Good analitycal properties will be necessary as we want to construct symbols of some Toeplitz operators out of elements of this group and then apply Widom's results (see below).
Given a matrix $g$ we denote with $\|g\|$ its Hilbert-Schmidt norm 
	$$\|g\|^2=\sum_{i,j=1}^{n}\|g_{i,j}\|^2$$

\begin{defin}
  Given a measurable matrix-valued loop $\gamma$ we define two norms $\|\gamma\|_\infty$ and $\|\gamma\|_{2,\unmez}$
	as
	$$\|\gamma\|_\infty:=\mathrm{ess}\sup_{\|z\|=1}\|\gamma(z)\|\quad\quad \|\gamma\|_{2,\unmez}:=\sum_k\Big(\left|k\right|\|\gamma^{(k)}\|^2\Big)^\unmez$$
	where we have Fourier expansion 
	$$\gamma(z)=\sum_{k=-\infty}^\infty{\gamma^{(k)}}z^k.$$
\end{defin}

\begin{defin}
$L_{\unmez}\mathrm{Gl}(n,\bC)$ is defined as the loop group of invertible measurable loops $\gamma$ such that
	$$\|\gamma\|_\infty+\|\gamma\|_{2,\unmez}<\infty.$$
\end{defin}

\begin{prop}[\cite{PS}]
 $L_\unmez \mathrm{Gl}(n,\bC)$ acts transitively on $\mathrm{Gr}^{(n)}$ and the isotropy group of $H_+^{(n)}$ is the group of constant loops $\mathrm{Gl}(n,\bC)$.
\end{prop}
 Proof can be found in \cite{PS}, here we just mention the principal steps necessary to arrive to this result.
 \begin{itemize}
  \item We define a subgroup $\mathrm{Gl}_{res}(H^{(n)})$ of invertible linear maps $g:H^{(n)}\rightarrow H^{(n)}$ acting on $\mathrm{Gr}(H^{(n)})$ (the restricted general linear group).
  \item We prove that every element of $\mathrm{Gl}_{res}(H^{(n)})$ commuting with multiplication by $z$ must belong to $L_\unmez \mathrm{Gl}(n,\bC)$.
  \item We take an element $W\in \mathrm{Gr}^{(n)}$ and a basis $\{w_1,...,w_n\}$ of the orthogonal complement of $zW$ in $W$.
  \item Out of this basis, putting vectors side by side, we construct $\cal W$ and easily check that $W=\cal W(z) H^{(n)}_+$.
  \item We verify that multiplication by $\cal W$ belongs to $\mathrm{Gl}_{res}(H^{(n)})$; since it obviously commutes with multiplication by $z$ we conclude that $\cal W(z)\in L_\unmez \mathrm{Gl}(n,\bC)$.
 \end{itemize}

\section{Szeg\"o-Widom theorem for block Toeplitz determinants.}

In his work (\cite{W0},\cite{W1} and \cite{W2}) H. Widom expressed the limit, for the size going to infinity, of certain block Toeplitz determinants as Fredholm determinants of an operator $\cal P$ acting on $H^{(n)}_+$.
Also he gave two different corollaries that allow us to compute this determinant in some particular cases.
In this section we recall, without proofs, these results.
Moreover we state Borodin-Okounkov formula as presented in \cite{BW} for matrix case.\\
We begin with some notations; given a loop $\gamma\in L_\unmez \mathrm{Gl}(n,\bC)$ we denote with $T_N(\gamma)$ the block Toeplitz matrix given by
$$
T_N(\gamma):=
\begin{pmatrix}
\gamma^{(0)} & \ldots & \ldots & \gamma^{(-N)}\\
&&&\\
\gamma^{(1)} & \ldots & \ldots & \gamma^{(-N+1)}\\
&&&\\
   \ldots       & \ldots & \ldots & \ldots\\
&&&\\
\gamma^{(N)} & \ldots & \ldots & \gamma^{(0)}
\end{pmatrix}$$
where we have the Fourier expansion $\gamma(z)=\sum_{k} \gamma^{(k)}z^k$.We denote $D_N(\gamma)$ its determinant.
We use the notation $T(\gamma)$ for the $\bN\times\bN$ matrix obtained letting $N$ go to infinity.
\begin{rem}
	It's easy to see that, in the base we have chosen above for $H^{(n)}$, $T(\gamma)$ is nothing but the matrix representation   of 
	$$\mathrm{pr}_+\circ\gamma:H^{(n)}_+\longrightarrow H^{(n)}_+$$
\end{rem}
\begin{thm}[Szeg\"o-Widom theorem,\cite{W2}]
		Suppose $\gamma\in L_\unmez \mathrm{Gl}(n,\bC)$ and  
		$$\underset{0\leq\theta\leq 2\pi}{\Delta}\arg\Big(\det\big(\gamma(e^{i\theta})\big)\Big)=0$$
		Then it exists the limit 
		$$D_\infty(\gamma):=\lim_{N\rightarrow\infty} \frac{D_N(\gamma)}{G(\gamma)^N}=\det(T(\gamma)T(\gamma^{-1}))$$
		where 
		$$G(\gamma)=\exp\Big(\frac{1}{2\pi}\int_0^{2\pi}\log\big(\det\gamma(e^{i\theta})\big)d\theta\Big)$$
\end{thm}

The proof of the theorem is contained in \cite{W2}; instead of rewriting it we simply consider the operator $T(\gamma)T(\gamma^{-1})$
and explain the meaning of "$\det$" in this case.
\begin{lemma}\label{ToeplitzHaenkel}
		Consider $\gamma _1,\gamma _2\in L_\unmez \mathrm{Gl}_n(n,\bC)$; we have
		$$T(\gamma _1\gamma _2)-T(\gamma _1)T(\gamma _2)=\Big[\sum_{k\geq 1}\gamma _1^{(i+k)}\gamma _2^{(-j-k)}\Big]_{i,j\geq 0}.$$
\end{lemma}
\prf 
The $(i,j)$-entry of left hand side reads
$$\sum_{k=-\infty}^\infty\gamma_1^{(i-k)}\gamma_2^{(k-j)}-
\sum_{k=0}^\infty\gamma_1^{(i-k)}\gamma_2^{(k-j)}=\sum_{k=\infty}^{-1}\gamma_1^{(i-k)}\gamma_2^{(k-j)}
=\sum_{k=0}^\infty\gamma_1^{(i+k+1)}\gamma_{2}^{(-k-j-1)}.$$ 
\endprf\\
In particular choosing $\gamma_1=\gamma$ and $\gamma_2=\gamma^{-1}$ we obtain
$$I-T(\gamma)T(\gamma^{-1})=\Big[\sum_{k\geq 1}\gamma^{(i+k)}(\gamma^{-1})^{(-j-k)}\Big]_{i,j\geq 0}$$
\begin{defin}
	\begin{equation}\label{P}
		\cal P_\gamma:=
		T(\gamma)T(\gamma^{-1})=\Bigg[\delta _i^j-\Big(\sum_{k\geq1}\gamma^{(i+k)}(\gamma^{-1})^{(-j-k)}\Big)\Bigg]_{i,j\geq 0}
	\end{equation}
\end{defin}
Thanks to the fact that 
$$\sum_{i\geq 0}\sum_{k\geq 1}\|\gamma^{(i+k)}\|^2=\sum_{k\geq 1}k\|\gamma^{(k)}\|^2<\infty$$
the product we have written on the right of (\ref{P}) is a product of two Hilbert-Schmidt operators.
So $\cal P_\gamma$ differs from the identity 
by a nuclear operator. Hence its determinant is well defined (see for instance \cite{Si}).
In our notation we obtained the equality
\begin{equation}\label{Widom}
	D_{\infty}(\gamma)=\det(\cal P_\gamma)
\end{equation}

We will call $\cal P_\gamma$ Plemelj's operator as it is related in a clear way with a Riemann-Hilbert factorization problem   (see section 4) already considered by Josip Plemelj in 1964 \cite{P}.\\
Unfortunately, in concrete cases, $\det(\cal P_\gamma)$ turns out to be really hard to compute; nevertheless we can use some shortcuts also provided by Widom in his works (\cite{W0},\cite{W1} and \cite{W2}).
\begin{prop}[\cite{W0}]\label{truncation}
	Suppose that $\gamma$ satisfies conditions imposed in Szeg\"o-Widom theorem and, moreover,
	$\gamma^{(i)}=0$ 	for $i\geq j+1$ or $\gamma^{(i)}=0$ for $i\leq j+1$.\\
	Then
	\begin{equation} \label{truncated Widom}
			D_\infty(\gamma)=D_j(\gamma^{-1})G(\gamma)^j
	\end{equation}
\end{prop}

\begin{prop}[\cite{W2}]\label{factorization}
 Suppose we have a symbol $\gamma$ satisfying conditions imposed in Szeg\"o-Widom theorem.
 Suppose moreover that $\gamma$ depends on a parameter $x$ in such a way that the function $x\rightarrow\gamma(x)$ is  differentiable.
 If $\gamma^{-1}$ admits two Riemann-Hilbert factorizations 
 $$\gamma^{-1}(z)=t_+(z)t_-(z)=s_-(z)s_+(z)$$
 such that
 \begin{gather*}
 	t_+(z):=\sum_{k\geq 0}t_+^{(k)}z^k\quad\quad s_+(z):=\sum_{k\geq 0}s_+^{(k)}z^k\\
 	t_-(z):=\sum_{k\leq 0}t_-^{(k)}z^k\quad\quad s_-(z):=\sum_{k\leq 0}s_-^{(k)}z^k
 \end{gather*}
 Then
 \begin{equation}\label{RH Widom}
		\frac{d}{dx}\log(D_\infty(\gamma))=
		\frac{i}{2\pi}\oint \mathrm{trace}\Bigg[\Big((\partial _zt_+)t_--(\partial _zs_-)s_+\Big)\partial _x\gamma\Bigg]dz.
	\end{equation}
\end{prop}

Also $D_N(\gamma)$ can be expressed as a Fredholm determinant as pointed out for the scalar case in \cite{BO} and generalized for matrix case in \cite{BW}.
\begin{thm}[Borodin-Okounkov formula, \cite{BW}]\label{BorodinOkounkov}
	Suppose that our symbol $\gamma(z)$ satisfying conditions of Szeg\"o-Widom's theorem admits two Riemann-Hilbert factorizations
	$$\gamma(z)=\gamma_+(z)\gamma_-(z)=\theta_-(z)\theta_+(z)$$
 such that
 \begin{gather*}
 	\gamma_+(z):=\sum_{k\geq 0}\gamma_+^{(k)}z^k\quad\quad \theta_+(z):=\sum_{k\geq 0}\theta_+^{(k)}z^k\\
 	\gamma_-(z):=\sum_{k\leq 0}\gamma_-^{(k)}z^k\quad\quad \theta_-(z):=\sum_{k\leq 0}\theta_-^{(k)}z^k
 \end{gather*}
 and $G(\gamma)=1$. Then for every $N$
 \begin{equation}
 	D_N(\gamma)=D_\infty(\gamma)\det(I-K_{\gamma,N})
 \end{equation}
 where, in coordinates, we have
 \begin{equation*}
 	(K_{\gamma,N})_{ij}=\begin{cases}
 												0\quad \text{if}\:\min\lbrace i,j\rbrace< N\\
 												\\
 												\sum_{k=1}^\infty (\gamma_-\theta_+^{-1})^{(i+k)}(\theta_-^{-1}\gamma_+)^{(-j-k)}\quad\text{otherwise.}
 											\end{cases}
 	\end{equation*}
\end{thm}

\begin{rem}
	One can easily verify that $\theta_-^{-1}\gamma_+$ is the inverse of $\gamma_-\theta_+^{-1}$ so that, again, we deal with operators of type $T(\phi)T(\phi^{-1})$ with $\phi=\gamma_-\theta_+^{-1}$. Also we want to point out that the assumption $G(\gamma)=1$ is not necessary. The formula for $G(\gamma)\neq 1$ is written in \cite{Bot}; since in our case we will always have $G(\gamma)=1$ we wrote the formula as it was given in \cite{BW}.
\end{rem}

\section{$\tau$ functions for constrained $\mathrm{KP}$ and Gelfand-Dickey\\ 
hierarchies as block Toeplitz determinants.}

It is well known that given a point $W\in \mathrm{Gr}$ its time evolution 
$$W(t):=\exp(t_1z+t_2z^2+t_3z^3+\ldots)W$$
is nothing but $\mathrm{KP}$ flow. Moreover points $W\in \mathrm{Gr}^{(n)}$ correspond to solutions of Gelfand-Dickey hierarchies (this is the celebrated Grassmannian formulation of $\mathrm{KP}$ hierarchy due to M.Sato, see for instance \cite{S}). In their paper \cite{SW} Segal and Wilson gave a formula for the corresponding $\tau$ function $\tau _W$ as determinant (Fredholm determinant) of the projection of $W(t)$ onto $H_+$. Here we take a slightly different approach that generalizes what has been done by Itzykson and Zuber in the study of Witten-Kontsevich $\tau$ function in \cite{IZ} (see also \cite{DF} and \cite{D2}).
This approach allows us to define not just $\tau_W$ but also a sequence of $\lbrace\tau_{W,N}\rbrace_{N>0}$ approximating $\tau_W$ and being themselves $\tau$ functions for some reductions of $\mathrm{KP}$.\\
Suppose we have an element $W\in \mathrm{Gr}^{(n)}$; thanks to results stated in Section 1 we can represent this element as
$$W=\begin{pmatrix} 	w_{11} & \ldots & \ldots & w_{n1} \\ 
					\ldots & \ldots & \ldots & \ldots\\
					\ldots & \ldots & \ldots &  \ldots\\
					w_{1n} & \ldots & \ldots & w_{nn}
					 \end{pmatrix}
					 H^{(n)}_+=\cal W(z) H^{(n)}_+$$
with $\cal W(z)\in L_{\unmez}\mathrm{Gl}(n,\bC)$.\\
Also we assume that the matrix $\cal W(z)=\lbrace w_{ij}(z)\rbrace_{i,j=1..n}$ satisfies
\begin{equation*}
		\begin{cases}
			w_{ii}=1+O(\frac{1}{z})\\
			\\
			w_{ij}=z(O(\frac{1}{z})),i>j\\
			\\
			w_{ij}=O(\frac{1}{z}),i<j
		\end{cases}
	\end{equation*}
This means that we restrict to the big cell, i.e. we assume $W\cong H^{(n)}_+$. In the sequel we will always assume that $W$ belongs to the big cell.
Obviously we have a base for $W\in \mathrm{Gr}^{(n)}$ given by 
$$\lbrace z^s w_j:s\in\bN,j=1\ldots n \rbrace$$ 
where $w_j$ is the column vector $(w_{1j}...w_{nj})^{\mathrm{T}}$. Using the isomorphim $\Xi:H^{(n)}\rightarrow H$ the corresponding base for $W\in \mathrm{Gr}$ is given by 
$$\lbrace \omega _{ns+j}=z^{ns} \Xi(w_j):s\in\bN,j=1\ldots n \rbrace$$ 
and, as in Section 1, we have
$$[\Xi(w_j)](z)=\sum_{i=1}^n z^{i-1}w_{ji}(z^n)$$
Note that thanks to the big-cell assumption we have $\omega _{ns+j}(z)=z^{ns+j-1}(1+O(\frac{1}{z}))$.\\
For these points $W\in \mathrm{Gr}^{(n)}$ and vectors spanning them we define the standard time evolution ($\mathrm{KP}$ flow) given by
$$\omega _{ns+j}(t;z):=\exp\Bigg(\sum_{i>0}t_iz^i\Bigg)\omega _{ns+j}(z)=\exp(\xi(t,z))\omega _{ns+j}(z)$$

Now we want to define the $\tau$ function associated to $W$ as limit for $N\rightarrow\infty$ of some block Toeplitz determinants $\tau_{W,N}$.
\begin{defin}\label{definitiontau}
	Take $M=Nn$ a multiple of  $n$.
	\begin{equation}\label{determinant}
			\tau_{W,N}(t):=\det\Bigg[\oint z^{-i}\omega _{j}(t;z)dz\Bigg]_{1\leq i,j\leq M=Nn}
	\end{equation}
\end{defin}

Fist of all we want to prove that $\tau_{W,N}$ is a block Toeplitz determinant and write explicitely the symbol.

\begin{lemma} For every $j=1\ldots n$ we have
	$$w_j(t;z):=\Xi^{-1}(\omega_j(t,z))=\exp(\xi(t,\Lambda))w_j(z)$$
	where we denote
		$$\Lambda:=\begin{pmatrix}
									0 & \ldots & \ldots & \ldots & z\\
									1 & 0   & \ddots & \ddots & 0\\
									0 & 1   & \ddots & \ddots & 0\\
									\ddots &\ddots& \ddots & \ddots & \ddots\\
									0 &\ddots& 0 & 1 & 0
								\end{pmatrix}$$
\end{lemma}
\prf
 We simply verify that multiplication by $z$ on $\mathrm{Gr}$ corresponds to multiplication by $\Lambda$ on $\mathrm{Gr}^{(n)}$ through the isomorphism $\Xi^{-1}$.
\endprf
\begin{prop}
	$\tau _{W,N}$ is the N-truncated $(n\times n)-$block Toeplitz determinant with symbol 
		\begin{equation}\label{GDSymbols}
			\cal W(t;z):=\exp(\xi(t,\Lambda))\cal W(z)
		\end{equation}
\end{prop}
\prf
	Take $i,j\leq n$ and $s,v\leq N$; the $(i+sn,j+vn)$-entry of the matrix in the right hand side of (\ref{determinant}) is given by
	\begin{gather*}
		\oint z^{-i-sn}\omega _{j+vn}(t;z)dz=\oint z^{-i-sn}z^{vn} \omega_j(t;z)dz=\\
		\oint z^{-i+(v-s)n}\sum_{k\in\bZ,l=1..n}w_{jl}(t)^{(k)}z^{nk+l-1}dz=w_{ji}(t)^{(s-v)}
	\end{gather*}
	so that the right hand side of (\ref{determinant}) is the transposed of the $N$-truncated $n\times n$ block Toeplitz matrix with symbol $\cal W(t;z)$.
\endprf\\
In the sequel of this paper we will call such symbols Gelfand-Dickey (GD) symbols.\\
Now generalizing what has been done by Itzykson and Zuber in \cite{IZ} we expand $\tau _{W,N}(t)$ in characters.
In this way, assigning degree $m$ to $t_m$, it's possible to state a certain property of stability for $\tau _{W,N}(t)$; namely every term of degree less or equal to $Q$ will be independent of $N$ for $Q\leq N$.
\begin{prop}
	$$\tau _{W,N}(t)=\sum _{l_1,...,l_{nN}\geq 0}\Big(\prod _{i}\omega _i^{(-l_i+i-1)}\Big)\chi _{l_1,...,l_{nN}}(X)$$
	where $X=\mathrm{diag}(x_1,...,x_{nN})$ is related to times $\lbrace t_i\rbrace$ through Miwa's parametrization
	$$t_k:=\mathrm{trace}\Bigg(\frac{X^k}{k}\Bigg)$$
	and 
	\begin{gather*}
		\chi _{l_1,...,l_{nN}}(X):=\dfrac{\det\begin{pmatrix} x_1^{l_1+nN-1} & x_1^{l_2+nN-2} &\ldots& x_1^{l_{nN}}\\
                                              \ldots			& \ldots			&\ldots&\ldots\\
											  x_{nN}^{l_1+nN-1} & x_{nN}^{l_2+nN-2} &\ldots& x_{nN}^{l_{nN}}
											  \end{pmatrix}}{\det\begin{pmatrix} x_1^{nN-1} & x_1^{nN-2} &\ldots& 1\\
                                              \ldots			& \ldots			&\ldots&\ldots\\
											  x_{nN}^{nN-1} & x_{nN}^{nN-2} &\ldots& 1
											  					\end{pmatrix}}
	\end{gather*}
\end{prop} 
\prf 
	We start from determinant representation (\ref{determinant}).
	The $(i,j)$-entry of the matrix will be 
	$$\sum_n\omega _j^{(-n)}p_{n+i-1}(t)$$ 
	where for every $n\geq 0$
	$$p_n(t):=\frac{1}{2\pi i}\oint \frac{\exp(\xi(t,z))}{z^{n+1}}dz$$\\
	are the classical Schur polynomials and $p_n(t)=0$ for every negative $n$.
	Then resumming everything we obtain
	$$\tau _{W,N}(t)=\sum _{k_1,\ldots,k_{nN}}\Big(\prod _j\omega _j^{(-k_j)}\Big)\det[p_{k_j+i-1}(t)]_{i,j=1\ldots nN}$$
	with $k_j\geq 1-j$.\\
	Equivalently we write
	\begin{gather*}
		\tau _{W,N}(t)=\sum _{l_1,\ldots ,l_{nN}\geq 0}\Big(\prod _j\omega _j^{(-l_j+j-1)}\Big)\det[p_{l_j-j+i}(t)]_{i,j=1\ldots nN}.
	\end{gather*}
	On the other hand it's well known that under Miwa's parametrization this last determinant can be written as 
	$\chi_{l_1,\ldots,l_{nN}}(X)$ (see for instance \cite{IZ},\cite{DF}); this completes the proof. 
\endprf\\
Now it's easy to see that assigning degree $1$ to every $x_i$ (which is equivalent to assigning degree $m$ to $t_m$) we obtain 
$$\mathrm{deg}(\chi _{l_1...l_{nN}})=\sum_{i=1}^Ml_i$$
From this easily verified property we obtain the following
\begin{lemma}
	Suppose
	$\mathrm{deg}(\chi _{l_1,...,l_{nN}})=Q\leq nN$. Then, if the character is different from zero, we have
	$$\chi _{l_1...l_{nN}}=\chi _{l_1,...,l_Q,0,...,0}.$$
\end{lemma}
\prf 
	Suppose $l_j\neq 0$, $j> Q$ and $l_i=0\text{ }\forall\text{ } i>j$.\\
	The $j^{th}$ column of the matrix $[p_{l_j-j+i}(t)]$ has positive subscripts $l_1+j-1,l_2+j-2,\ldots ,l_j$.\\
	On the other hand $\sum l_i=Q$; hence the sum of these subscripts is 
	$$Q+\sum_{r=0}^{j-1}r\leq\sum_{r=0}^{j}r$$
	hence two subscripts must be equal, then two lines of the matrix are equal. 
\endprf\\
From this corollary it follows directly the following result.
\begin{prop}\label{stabilization}
	Up to degree $Q$ the function $\tau _{W,N}(t)$ does not depend on $N$ with $N\geq Q$.
\end{prop}
This proposition allows us to define in a rigorous sense
\begin{equation}\label{tau}
	\tau_W(t):=\lim_{N\rightarrow\infty}\tau_{W,N}(t)
\end{equation}
where the limit is defined as limit of formal graded series in $t_i$; this means that 
$$\lim_{N\rightarrow\infty}\deg(\tau_{W,N}(t)-\tau_{W}(t))=\infty$$
We will say that $\tau _W(t)$ is the \emph{stable limit} of $\tau _{W,N}(t)$ for $N\rightarrow\infty$.\\
On the other hand, in the sequel, we will prove that the symbol $\cal W(t)$ satisfies Szeg\"o-Widom's condition for every values of $t_i$ so that the limit in (\ref{tau}) exist pointwise in time parameters and can be written as a Fredholm determinant.\\
Now, following again \cite{IZ}, we write a differential operator $\Delta _{W,N}(t)$ associated to the function $\tau _{W,N}(t)$. In the sequel we will always write $D$ for the partial derivative with respect to $t_1$. We will prove that for every $N$ the pseudo-differential operator $\Delta _{W,N}(t)D^{-nN}$ satisfies Sato's equations for the dressing (see \cite{S}) and we recover the usual relation between $\tau$ and wave functions.
\begin{lemma}
	Define 
	$$f_{s,N}(t):=\sum_{k>s}\omega _s^{(-k)}p_{k+nN-1}(t)$$
	Then we have
	$$\tau _{W,N}(t)=\mathrm{Wr}(f_{1,N}(t),\ldots,f_{nN,N}(t)):=\det[D^{nN-j}f_{i,N}(t)]_{1\leq i,j\leq nN}$$
\end{lemma}
\prf
	From definition \ref{definitiontau} the $(i,j)$-entry of matrix defining $\tau_{W,N}(t)$ is
	$$\sum_{k>j}\omega _j^{(-k)}p_{k+i-1}(t)$$
	On the other hand we have
	\begin{gather*}
		D^{nN-j}f_{i,N}(t)=D^{nN-j}\Bigg(\sum_{k>i}\omega _i^{(-k)}p_{k+nN-1}(t)\Bigg)=\sum_{k>i}\omega _i^{(-k)}p_{k+j-1}(t)
	\end{gather*}
	(using the equation $D^s(p_m(t))=p_{m-s}(t)$). Hence we obtained the proof.
\endprf
\begin{defin}
	We define the differential operator $\Delta _{W,N}(t)$ of order $N$ in $D$ as
	$$\Delta _{W,N}(f):=\dfrac{\mathrm{Wr}(f,f_{1,N}(t),\ldots,f_{nN,N}(t))}{\mathrm{Wr}(f_{1,N}(t),\ldots,f_{nN,N}(t))}$$
	where $f\in H^{(n)}$ depends in a differentiable way on $\lbrace t_i\rbrace_{i\geq 1}$.
\end{defin}

\begin{prop}
	The following equations for time-derivatves of $\Delta_{W,N}(t)$ holds:
	\begin{equation}\label{Sato1}
		\dfrac{\partial}{\partial_{t_i}}(\Delta _{W,N}(t))=(\Delta _{W,N}(t)D^i\Delta _{W,N}^{-1}(t))_+-\Delta _{W,N}(t)D^i
	\end{equation}
\end{prop}
\prf
	It is enough to prove the equality of the two differential operators when acting on $f_{1,N}(t),...f_{nN,N}(t)$ which are $nN$ indipendent 	solutions of the equation
	$$(\Delta _{W,N}(t))(f(t))=0$$
	But this amounts to proving
	$$\Bigg[\dfrac{\partial}{\partial_{t_i}}(\Delta _{W,N})(t)\Bigg]f_{j,N}(t)+\Delta _{W,N}(t)\dfrac{\partial^i}{\partial t_1^i}f_{j,N}(t)=0\quad\forall j$$
	which is true iff
	$$\dfrac{\partial}{\partial t_i}(\Delta _{W,N}(t)f_{j,N}(t))=0\quad\forall j.$$
	This equality is obviously satisfied. 
\endprf\\
Multiplying $\Delta_{W,N}$ from the right with $D^{-nN}$ we found a pseudodifferential operator that, in fact, gives a solution of $\mathrm{KP}$ equations.
\begin{defin}
	$$S_{W,N}(t):=\Delta_{W,N}(t)D^{-nN}$$
\end{defin}	
\begin{prop}
	$S_{W,N}(t)$ is a monic pseudo-differential operator of order $0$ satisying Sato's equation (see for instance \cite{S})
	\begin{equation}\label{Sato2}
		\dfrac{\partial}{\partial_{t_i}}(S_{W,N}(t))=(S_{W,N}(t)D^iS_{W,N}^{-1}(t))_+-S _{W,N}(t)D^i
	\end{equation}
	Hence the monic pseudo-differential operator of order $1$
	\begin{equation}\label{pseudodifferential}
		\cal L_{W,N}(t):=(S_{W,N}(t)DS_{W,N}^{-1}(t))
	\end{equation}
	satisfies the usual Lax system for $\mathrm{KP}$
	\begin{equation}\label{Lax}
		\frac{\partial\cal L_{W,N}}{\partial t_k}=[(\cal L^k_{W,N})_+,\cal L_{W,N}]
	\end{equation}
\end{prop}
\prf
	It is obvious that $S_{W,N}$ is a monic pseudo-differential operator of order $0$ since $\Delta _{W,N}$, which is of order $nN$, is normalized so that the leading term is equal to $1$.
	Equation (\ref{Sato2}) follows directly from (\ref{Sato1}).
	The derivation of Lax system from Sato's equations is well known: one has just to derive the relation
	$$\cal L_{W,N}S_{W,N}=S_{W,N}D$$
	for $t_k$ and use the obvious relation $[\cal L_{W,N},\cal L_{W,N}^k]=0$
\endprf\\
It remains to prove that $\tau_{W,N}(t)$ is really the $\tau$ function for these solutions $\cal L_{W,N}(t)$ of $\mathrm{KP}$ equations.
We recall the usual relations between the dressing $S$, the wave function $\psi$ and $\tau$ function given by
$$\psi(t;z)=S(t)(\exp(\xi(t,z)))=\exp(\xi(t,z))\frac{\tau(t-\frac{1}{[z]})}{\tau(t)}$$
(we recall that the notation $t-\frac{1}{[z]}$ stands for the vector with $i^{th}$ component equal to $t_i-\frac{1}{i z^i}$)
All we have to prove is the following

\begin{prop}
	\begin{equation}\label{wavefunction}
		\psi_{W,N}(t;z):=S_{W,N}(t)(\exp(\xi(t,z)))=\exp(\xi(t,z))\frac{\tau _{W,N}(t-\frac{1}{[z]})}{\tau_{W,N}(t)}
	\end{equation}
\end{prop}
\prf
	Equivalently we prove that
	$$(\Delta _{W,N}(t))\exp(\xi(t,z))=\exp(\xi(t,z))z^{nN}\dfrac{\tau _{W,N}(t-\frac{1}{[z]})}{\tau _{W,N}(t)}$$
	Since we have 
	$$p_n(t-\frac{1}{[z]})=p_n(t)-z^{-1}p_{n-1}(t)$$ 
	the right hand side of the equality above can be written as
	$$z^{nN}e^{\xi(x,t)}\frac{\det\begin{pmatrix}
			D^{nN-1}f_1-z^{-1}D^{nN}f_1&\ldots& f_1-z^{-1}Df_1\\
			\ldots                        &\ldots&\ldots\\
			D^{nN-1}f_{nN}-z^{-1}D^{nN}f_{nN}&\ldots& f_{nN}-z^{-1}Df_{nN}
             \end{pmatrix}}{\mathrm{Wr}(f_1,\ldots,f_{nN})}$$
	(here derivative is with respect to $t_1$, we don't write dependence on $f_i$ on $t$ to avoid heavy notation)
	The left hand side can be written as 
	$$\frac{\det\begin{pmatrix}
		z^{nN}e^{\xi(t,z)}&\ldots&\ldots&e^{\xi(t,z)}\\
		D^{nN}f_1&\ldots&\ldots&f_1\\
		\ldots      &\ldots&\ldots&\ldots\\
		D^{nN}f_{nN}&\ldots&\ldots&f_{nN}
	\end{pmatrix}}{\mathrm{Wr}(f_1,\ldots,f_{nN})}.$$
	It is easy to check that these two expressions are equal.
\endprf\\
We want now to study the structure of $\cal L_{W,N}$ with more attention; our investigation will lead us to discover that, actually, we are dealing with rational reductions (\cite{Kr},\cite{D2}) of $\mathrm{KP}$.\\
First of all we recall a useful lemma (proof can be found for instance in \cite{I}).
\begin{lemma}\label{wronsky}
	Let $\lbrace g_1,...,g_m\rbrace$ be a basis of linearly independent solutions of a differential operator $K$ of order $m$.
	Then one can factorize $K$ as
	$$K=(D+T_m)(D+T_{m-1})\cdots(D+T_1)$$
	with
	$$T_j=\frac{\mathrm{Wr}(g_1,\ldots,g_{m-1})}{\mathrm{Wr}(g_1,\ldots,g_m)}.$$
\end{lemma}
We will state now properties of symmetry for $f_{s,N}$ that will be useful in the sequel.
\begin{prop}
	The following equalities hold:
		\begin{equation}\label{symmetry1}
			f_{s,N+1}(t)=f_{s-n,N}(t)
		\end{equation}
		\begin{equation}\label{symmetry2}
			f_{s,N}(t)=D^{n}f_{s+n,N}(t)
		\end{equation}
\end{prop}
\prf
	We will use the equality
	$$\omega_s^{(-k)}=\omega_{s+n}^{(-k+n)}$$
	which follows from the very definition of these coefficients.\\
	Then for (\ref{symmetry1}) we have
	\begin{gather*}
		f_{s,N+1}(t)=\sum_k\omega_s^{(-k)}p_{k+n+nN-1}(t)=\sum_k\omega_{s-n}^{(-k-n)}p_{k+n+nN-1}(t)=\\
		\sum_k\omega_{s-n}^{(-k)}p_{k+nN-1}(t)=f_{s-n,N}(t)
	\end{gather*}
	For (\ref{symmetry2}) we have
	\begin{gather*}
		D^n(f_{s+n,N}(t))=D^n\sum_k\omega_{s+n}^{(-k)}p_{k+nN-1}(t)=\sum_k\omega_{s+n}^{(-k)}p_{k+nN-1-n}(t)=\\
		\sum_{k}\omega_s^{(k-n)}p_{k-n+nN-1}(t)=f_{s,N}(t)
	\end{gather*}
\endprf

\begin{thm}\label{constrained}
	For every $N>0$ the pseudodifferential operator $\cal L_{W,N}$ and its $n^{th}$ power
	$$L_{W,N}=\cal L_{W,N}^n$$
	can be factorized as
	\begin{itemize}
		\item $\cal L_{W,N}=\cal L_{1,W,N}(\cal L_{2,W,N})^{-1}$
		\item $L_{W,N}=M_{1,W,N}(M_{2,W,N})^{-1}$
	\end{itemize}
	where all the factors are differential operators and 
	\begin{itemize}
		\item $\mathrm{ord}(\cal L_{1,W,N})=nN+1, \quad\mathrm{ord}(\cal L_{2,W,N})=nN$
		\item $\mathrm{ord}(M_{1,W,N})=2n, \quad\mathrm{ord}(M_{2,W,N})=n$
	\end{itemize}
\end{thm}
\prf
	The first factorization comes directly from the fact that
	$$\cal L_{W,N}(t)=\Delta_{W,N}(t)D(\Delta_{W,N}(t))^{-1}$$
	For the second factorization we note that we have the factorization
	$$L_{W,N}(t)=\Delta_{W,N}(t)D^n(\Delta_{W,N}(t))^{-1}$$
	where the first operator $\Delta_{W,N}(t)D^n$ has order $M+n$ while the second (i.e. $\Delta_{W,N}$) has order $nN$.
	Moreover as follows from (\ref{symmetry2}) we have
	\begin{itemize}
		\item $\Delta_{W,N}f_{i,N}=0\quad\forall i=1,\ldots,nN$
		\item $\Delta_{W,N}D^nf_{i,N}=0\quad\forall i=n+1,\ldots,nN$
	\end{itemize}
	hence using lemma \ref{wronsky} one can simplify factorization above as
	$$L_{W,N}=M_{1,W,N}(M_{2,W,N})^{-1}$$
	where $M_{2,W,N}$ is given explicitely by the formula
	$$M_{2,W,N}=(D+K_{n,N})(D+K_{n-1,N})...(D+K_{1,N})$$
	with
	$$K_{j,N}=
	D\Bigg[ \log\Bigg(\dfrac{\mathrm{Wr}(f_{n+1,N},\ldots,f_{nN,N},f_{1,N},\ldots,f_{j-1,N})}{\mathrm{Wr}(f_{n+1,N},\ldots,f_{nN,N},f_{1,N},\ldots,f_{j,N})}\Bigg)\Bigg]$$  
\endprf
\begin{thm}\label{rationalreduction}
 The sequence $\lbrace\cal L_{W,N}\rbrace_{N\geq 1}$ satisfies recursion relation
 \begin{equation}\label{recursion}
 	\cal L_{W,N+1}=\cal T_N\cal L_{W,N}(\cal T_N)^{-1}
 \end{equation}
 with
 $$\cal T_N=(D+T_{n,N})(D+T_{n-1,N})...(D+T_{1,N})$$
 $$T_j=
 D \log\Bigg(\dfrac{\mathrm{Wr}(f_{1,N},...,f_{nN,N},f_{1,N+1},...f_{j-1,N+1})}{\mathrm{Wr}(f_{1,N},...,f_{nN,N},f_{1,N+1},...f_{j,N+1})}\Bigg).$$
\end{thm}
\prf
	We observe that thanks to (\ref{symmetry1})
	\begin{itemize}
		\item $\Delta_{W,N}f_{i,N}=0\quad\forall i=1,\ldots,nN$
		\item $\Delta_{W,N+1}f_{i,N}=0\quad\forall i=n+1,\ldots,nN$
		\item $\Delta_{W,N+1}f_{i,N+1}=0\quad\forall i=1,\ldots,n$
	\end{itemize}
	Hence using again (\ref{wronsky}) we obtain the recursion relation
	$$\Delta_{W,N+1}=\cal T_N\Delta_{W,N}$$
	and from this last equation we recover the recursion relation for the Lax operator.
\endprf\\
Using the notation of \cite{Ar} we say that, for every $N$, $\tau_{W,N}\in\mathrm{cKP}_{1,nN}\cap\mathrm{cKP}_{n,n}$.\\
This means that given a $\tau$ function for $\mathrm{KP}$ with corresponding Lax pseudodifferential operator $\cal L$ we say that $\tau\in\mathrm{cKP}_{m,n}$ iff $\cal L^m$ can be written as the ratio of two differential operators of order $m+n$ and $n$ respectively. These special reductions of $\mathrm{KP}$ begun to be studied in 1995 by Dickey and Krichever (\cite{D2},\cite{Kr}); a geometric interpretation of corresponding points in the Grassmannian has been given in \cite{vdL1} and \cite{vdL2}.
Going back to theorem \ref{rationalreduction} we point out that the first decomposition as well as the recursion formula are already known and, as pointed out in \cite{Ar}, come simply from the fact that we have a truncated dressing.
Actually our sequence of $\lbrace\tau _{W,N}\rbrace_{N\geq 1}$ is a part of a sequence already studied by Dickey in \cite{D1} under the name of stabilizing chain; in that article Dickey already provided the recursion formula written above as well as some differential equations for coefficients of $\cal T_N$.
Nevertheless, to our best knowledge, connection with block Toeplitz determinants never appeared before.
Also the fact that $\tau_{W,N}\in\mathrm{cKP_{n,n}}$ is something new.
It could be interesting to find recursion relations as well as differential equations for $M_{1,W,N}$ and $M_{2,W,N}$; we plan to do it in a subsequent work.
Till now all we can do is to infer from recursion formula for Lax operator the following formula
\begin{equation}
	(M_{1,W,N+1})^{-1}\cal T_{N}M_{1,W,N}=(M_{2,W,N+1})^{-1}\cal T_{N}M_{2,W,N}.
\end{equation}
Now we want to go one step further and see what happens for $N\rightarrow\infty$.
Obviously thanks to the property of stabilization stated in proposition \ref{stabilization} we can define a pseudodifferential operator $\cal L_W$ and a wave function $\psi _W$ related to $\tau _W$ in the same way as for finite $N$ and we will obtain a solution of $\mathrm{KP}$ as well.
Actually a stronger statement holds.
\begin{prop}\label{nKdV}
	Given $W\in \mathrm{Gr}^{(n)}$ the functions $\tau _W$, $\psi _W$ and $L_W:=(\cal L_W)^n$ are respectively the $\tau$ function, the wave function and the differential operator of order $n$ corresponding to a solution of $n^{th}$ Gelfand-Dickey hierarchy.
\end{prop}
\prf
	It is known \cite{SW} that subspaces satisfying $z^nW\subseteq W$ correspond to solutions of $n^{th}$ Gelfand-Dickey hierarchy.
	What we have to prove is that $L_W(t)=(L_W(t))_+$.\\
	In the sequel we will omit dependence on times of $L_W$ and $S_W$ if no ambiguity arises.\\
	From the usual relation
	\begin{equation}\label{wave}
		\frac{\partial\psi_W}{\partial t_n}(t;z)=(L_W(t))_+\psi_W(t;z)
	\end{equation}
	we obtain immediately
	$$\frac{\partial S_W}{\partial t_n}+S_WD^n=(L_W)_+S_W$$
	so that we have to prove that 
	$$\frac{\partial S_W}{\partial t_n}=0$$
	On the other hand
	$$\psi_W(t;z)=\exp(\xi(t,z))(1+\sum_{1=1}^\infty s_i(t)z^{-i})$$
	where
	$$S_W(t)=1+\sum_{1=1}^\infty s_i(t)D^{-i}$$
	Using this explicit expression for the wave function and substituting in (\ref{wave}) we obtain
	$$(L_W(t))_+\psi_W(t;z)-z^n\psi_W(t;z)=\exp{(\xi(t,z))}\sum_{i=1}^{\infty}\frac{\partial s_i(t)}{\partial t_n}z^{-i}$$
	The left hand side of this equation lies on $W(t)=\exp(\xi(t,z))W$ for every $t$ so that multiplying both terms for $\exp({-\xi(t,z)})$
	one obtains that they belong to subspaces transverse one to the other ($W$ and $H_-$), hence both of them vanish. This means that $\frac{\partial s_i}{\partial t_n}=0$ for every $i$.
\endprf\\
In virtue of this proposition, when computing $\tau _W$ associated to $W\in \mathrm{Gr}^{(n)}$, we will always omit times $t_{jn}$ multiple of $n$. Setting $\lbrace t_{jn}=0,j\in\bN\rbrace$ will be important in order to be able to apply Szeg\"o-Widom's theorem; in this case we will write $\tilde t$ instead of $t$.
\begin{prop}
	Take any $W\in \mathrm{Gr}^{(n)}$ in the big cell of $\mathrm{Gr}^{(n)}$ and a corresponding GD symbol $\cal W(t;z)$.\\
	Then
	\begin{equation}\label{main}
	\tau_{W}(\tilde t)=\det(\cal P_{\cal W(\tilde t;z)}).
	\end{equation}
\end{prop}
\prf
	All we have to prove is that conditions of Szeg\"o-Widom's theorem are satisfied and $G(\cal W(\tilde t;z))=1$.
	We observe that 
	$$\cal W(\tilde t;z)\in L_\unmez \mathrm{Gl}(n,\bC)\quad \forall \tilde t$$
	since we can always find $\cal W(z)\in L_\unmez \mathrm{Gl}(n,\bC)$ such that $W=\cal W(z)H^{(n)}_+$ and $\exp(\xi(\tilde t,\Lambda))$ is continuously differentiable (obviously when restricted to a finite number of times).\\
	Moreover 
	$$\det[\exp(\xi(\tilde t,\Lambda))]=1$$
	 since we deleted times multiple of $n$ and 
	 $\det(\cal W(z))=1+O(z^{-1})$ by big cell assumption.\\
	This implies that we have 
	$$\underset{0\leq\theta\leq 2\pi}{\Delta}\det\Big(\cal W(\tilde t;e^{i\theta})\Big)=0$$ and 
	$$G(\cal W(\tilde t;z))=1.$$ 
\endprf\\
We are now in the position to state the main result of this paper, proof follows from results obtained above.
\begin{thm}
	Given any point $W\in \mathrm{Gr}^{(n)}$ and corresponding GD symbol $\cal W(t;z)$ the following facts hold true:
	\begin{itemize}
		\item $\lbrace\tau _{W,N}(t):=D_N(\cal W(t;z))\rbrace_{0\leq N<\infty}$ \\
		is a sequence of $\tau$ functions for $\mathrm{KP}$ associated to wave function $\psi_{W,N}(t;z)$ and \\
		pseudodifferential operators $\cal L_{W,N}(t)$ given respectively by (\ref{wavefunction}) and (\ref{pseudodifferential}).
		\item For every $N>0$ we have $\tau_{W,N}\in \mathrm{cKP}_{1,nN}\cap \mathrm{cKP}_{n,n}$.\\
		 Explicit factorizations of Lax operator and its $n^{th}$ power are given in theorem \ref{constrained}.
		\item The sequence admits stable limit $\tau _W(\tilde t)$.\\
		The stable limit is a solution of $n^{th}$ Gelfand-Dickey hierarchy. It is equal to the Fredholm determinant 
		$\det{\cal P_{\cal W(\tilde t;z)}}$.
	\end{itemize}
\end{thm}

\begin{rem}
Also all the $\tau_{W,N}(\tilde t)$ can be expressed as Fredholm determinants; in order to give explicit expressions we need  a certain Riemann-Hilbert factorization of symbol $\cal W(\tilde t;z)$. This factorization will be obtained in section 5 and it will be exploited to express $\tau_{W,N}(\tilde t)$ as a Fredholm determinant.
\end{rem}

\section{Riemann-Hilbert problem and Plemelj's integral formula.}

It is evident from proposition \ref{factorization} that Riemann-Hilbert decompositions of symbol $\gamma$ for a block Toeplitz operator plays an important role in computing $D_\infty(\gamma)$.\\
Here we will show that actually Plemelj's operator itself enters in a integral equation (see \cite{P}) giving solutions of Riemann-Hilbert problem
\begin{equation}\label{Plemelj RH}
	\varphi_+(z)=\gamma^{\mathrm{T}}(z)\varphi_-(z).
\end{equation}
Here $\varphi_+(z)$ and $\varphi_-(z)$ are respectively analytical functions defined inside and outside the circle.
In this section we consider a smaller class of loops; $\gamma(z)$ will be a matrix-valued function that extends analytically on a neighborhood of $S^1$. For convenience of the reader we recall here main steps to arrive to Plemelj's integral formula \cite{P}.
\begin{lemma}
	Suppose that $f_+(z),f_-(z)$ are functions on $S^1$ satisfying 
$$\left|f(\zeta_2)-f(\zeta_1)\right|<\left|\zeta_2-\zeta_1\right|^\mu C$$
for some positive constants $\mu,C$ and for every $\zeta_1,\zeta_2\in S^1$.
	Necessary and sufficient conditions for $f_+(z)$ and $f_-(z)$ to be boundary values of analytic functions regular inside or outside $S^1\subseteq\bC$ and with value $c$ at infinity are respectively
	\begin{align}
		\frac{1}{2\pi i}\oint\frac{f_+(\zeta)-f_+(z)}{\zeta-z}d\zeta=0\label{Cauchy1}\\
		\frac{1}{2\pi i}\oint\frac{f_-(\zeta)-f_-(z)}{\zeta-z}d\zeta+f_-(z)-c=0\label{Cauchy2}
	\end{align}
\end{lemma}
We have to point out that here both $\zeta$ and $z$ lies on $S^1$ so that one has to be careful and define (\ref{Cauchy1}) and (\ref{Cauchy2}) as appropriate limits.
Namely one proves that taking $\zeta$ slightly inside or outside $S^1$ along the normal and making it approach to the circle we obtain the same result which will be, by definition, the value of our integral.
Now suppose we want to find solutions of (\ref{Plemelj RH}); we normalize the problem requiring $\varphi_-$ taking value $C$ at infinity.
Taking an appropriate linear combination of (\ref{Cauchy1}) and (\ref{Cauchy2}) and using (\ref{Plemelj RH}) we find that $\varphi_-(z)$ must satisfy the equation
\begin{equation}\label{integralPlemelj1}
	C=
	\varphi(z)-\frac{1}{2\pi i}\oint\frac{(\gamma^{\mathrm{T}})^{-1}(z)\gamma^{\mathrm{T}}(\zeta)-I}{\zeta-z}\varphi(\zeta)d\zeta
\end{equation}

Note that here we do not have to take any limit since the integrand is well defined for every point of $S^1$.
We also want to consider the associate homogeneous equation
\begin{equation}\label{integralPlemelj2}
	0=
	\varphi(z)-\frac{1}{2\pi i}\oint\frac{(\gamma^{\mathrm{T}})^{-1}(z)\gamma^{\mathrm{T}}(\zeta)-I}{\zeta-z}\varphi(\zeta)d\zeta
\end{equation} 
as well as its adjoint
\begin{equation}\label{integralPlemelj3}
	0=\psi(z)+\frac{1}{2\pi i}\oint\frac{\gamma(z)\gamma^{-1}(\zeta)-I}{\zeta-z}\psi(\zeta)d\zeta
\end{equation}
Obviously, as usual in Fredholm's theory, the equations (\ref{integralPlemelj2}) and (\ref{integralPlemelj3}) either have only  trivial solution or they have the same number of linearly independent solutions.
\begin{lemma}
	Consider two adjoint RH problems
	\begin{align}
		\varphi_+(z)=\gamma(z)^{\mathrm{T}}\varphi_-(z)\label{hom1}\\
		\psi_+(z)=\gamma(z)\psi_-(z)\label{hom2}
	\end{align}
	normalized as $\psi_-(\infty)=\varphi_-(\infty)=0$.\\
	Any solution $\varphi_-$ of (\ref{hom1}) is a solution of (\ref{integralPlemelj2}) as well as any solution $\psi_+$ of (\ref{hom2}) is a solution of (\ref{integralPlemelj3}).  
\end{lemma}
\prf 
	We just repeat computations made for non-homogeneous case.
\endprf\\
Now we introduce a new integrable operator acting on $H^{(n)}_+$ and prove that it is actually equal to the Plemelj's operator.
\begin{defin}
	For every $f\in H^{(n)}_+$ we define
	\begin{equation}\label{newPlemeljoperator}
		[\tilde {\cal P}_\gamma \psi](z):=
		\mathrm{pr}_+\Bigg(\psi(z)+\frac{1}{2\pi i}\oint\frac{\gamma(z)\gamma^{-1}(\zeta)-I}{\zeta-z}\psi(\zeta)d\zeta\Bigg)
	\end{equation}
	where $\mathrm{pr}_+$ denote the projection onto $H^{(n)}_+$.
\end{defin}

\begin{prop}
	$$\tilde {\cal P}_\gamma=\cal P_\gamma.$$
\end{prop}

\prf 
	We write $\tilde {\cal P}_\gamma$ in coordinates and verify we obtain the same as in (\ref{P}).
	To do so as in the definition of integrals (\ref{Cauchy1}) and (\ref{Cauchy2}) we compute (\ref{newPlemeljoperator}) imposing $\left|\zeta\right|<\left| z \right|$; the formula will hold when $\zeta$ approach to $S^1$ in the same way as in (\ref{Cauchy1}) and (\ref{Cauchy2}).
	For a consistency check we will prove we obtain the same result imposing $\left|\zeta\right|>\left| z \right|$.
	Let's start with $\left|\zeta\right|<\left| z \right|$; we have
	\begin{gather*}
		\psi(z)+\frac{1}{2\pi i}\oint\frac{\gamma(z)\gamma^{-1}(\zeta)-I}{\zeta-z}\psi(\zeta)d\zeta=\\
		\\
		\psi(z)+\frac{1}{2\pi i}
		\oint\sum_{k\geq 1}\frac{\zeta^k}{z^k}\Big(I-\sum_{p,q\in\bZ}\gamma^{(p)}(\gamma^{-1})^{(q)}z^p\zeta^q\Big)
		\sum_{s\geq 0}\psi^{(s)}\zeta^s\frac{d\zeta}{\zeta}
	\end{gather*}
	Imposing $k+q+s=0$ we get that this is equal to
	\begin{gather*}
		\psi(z)+\sum_{p\in\bZ}\sum_{k\geq 1}\sum_{s\geq 0}\gamma^{(p)}(\gamma^{-1})^{(-k-s)}\psi^{(s)}z^{p-k}=
		\psi(z)+\sum_{t\in\bZ}\sum_{k\geq 1}\sum_{s\geq 0}\gamma^{(t+k)}(\gamma^{-1})^{(-k-s)}\psi^{(s)}z^t
	\end{gather*}
	Taking the projection on $H_+^{(n)}$ we obtain exactly formula (\ref{P}).
	Now for $\left|\zeta\right|>\left| z \right|$ we have
	\begin{gather*}
		\psi(z)+\frac{1}{2\pi i}\oint\frac{\gamma(z)\gamma^{-1}(\zeta)-I}{\zeta-z}\psi(\zeta)d\zeta=\\
		\\
		\psi(z)+\frac{1}{2\pi i}
		\oint\sum_{k\geq 0}\frac{\zeta^k}{z^k}\Big(\sum_{p,q\in\bZ}\gamma^{(p)}(\gamma^{-1})^{(q)}z^p\zeta^q-I\Big)
		\sum_{s\geq 0}\psi^{(s)}\zeta^s\frac{d\zeta}{\zeta}
	\end{gather*}
	Imposing $q+s=k$ we arrive to
	\begin{gather*}
		\sum_{k,s\geq 0}\sum_{p\in\bZ}\gamma^{(p)}(\gamma^{-1})^{(k-s)}\psi^{(s)}z^{k+p}=
		\sum_{k,s\geq 0}\sum_{t\in\bZ}\gamma^{(t-k)}(\gamma^{-1})^{(k-s)}\psi^{(s)}z^{t}
	\end{gather*}
	Taking the projection on $H_+^{(n)}$ we obtain that this is equal to $T(\gamma)T(\gamma^{-1})$ so that the two computations for $\left|\zeta\right|<\left| z \right|$ and for $\left|\zeta\right|>\left| z \right|$ coincide in virtue of lemma \ref{ToeplitzHaenkel}
\endprf

\begin{thm}
	Suppose we are given a symbol $\gamma(z)$ analytic in a neighborhood of $S^1$ and such that 
	$$D_\infty(\gamma)\neq 0$$
	Then the Riemann-Hilbert problem 
	$$\varphi_+(z)=\gamma(z)^{\mathrm{T}}\varphi_-(z)$$ 
	normalized as $\varphi_-(\infty)=C$ 
	admits (if existing) a unique solution. 
\end{thm}
\prf
	Suppose we have two distinct solutions $(\varphi_{1-},\varphi_{1+})$ and $(\varphi_{2-},\varphi_{2+})$; taking the difference we obtain a non-trivial solution of (\ref{hom1}).
	Then also (\ref{hom2}) admits non trivial solutions and the same holds for (\ref{integralPlemelj3}).
	But this means that we have a non zero $\psi(z)\in H^{(n)}_+$ such that $[\cal P_\gamma\psi](z)=0$ which is impossible since
	$$\det(\cal P_\gamma)=D_\infty(\gamma)\neq 0$$
\endprf\\
Existence of factorization will be treated in the next section for the specific case of Gelfand-Dickey symbols.
For a general treatment of the problem of existence see \cite{P}.

\section{Factorization for Gelfand-Dickey symbols}

Here we will prove that for Gelfand-Dickey symbols we can write the unique solution of factorization (\ref{Plemelj RH}) in terms of data
$L_W(\tilde t),\psi_W(\tilde t;z)$. We recall that $L_W(\tilde t)$ and $\psi_W(\tilde t;z)$ are the stable limits of $L_{W,N}(\tilde t)$ and $\psi_{W,N}(\tilde t;z)$. They represent the differential operator and the wave function associated to the solution $\tau_W(\tilde t)$.
Our exposition here is closely related to \cite{SS}. At the end of the section we will use the factorization obtained to express any $\tau_{W,N}(\tilde t)$ as a Fredholm determinant.
As we have written before in the proof of proposition \ref{nKdV} we have the relation
\begin{equation}\label{eigenequation}
	L_W(\tilde t)\psi_W(\tilde t;z)=z^n\psi_W(\tilde t;z)
\end{equation}
where $\psi_W(\tilde t;z)$ admits asymptotic expansion
$$\psi_W(\tilde t;z)=\exp(\xi(\tilde t,z))(1+O(z^{-1}))$$
Now out of $\psi_W$ we construct $n$ time-dependent functions
$$\psi_{W,i}(\tilde t;z):=D^i(\psi_W(\tilde t;z)):i=0,\ldots ,n-1$$
belonging to the subspace $W\in \mathrm{Gr}$.
\begin{defin}
	$$\Psi_W(\tilde t;z):=\begin{pmatrix}
					1 & \zeta_1 & \ldots & \zeta_1^{n-1}\\
					&&&\\
					1 & \zeta_2 & \ldots & \zeta_2^{n-1}\\
					&&&\\
					\ldots &\ldots & \ldots & \ldots\\
					&&&\\
					1 & \zeta_n & \ldots & \zeta_n^{n-1}
			\end{pmatrix}^{-1}
			\begin{pmatrix}
					\psi_{W,0}(\tilde t;\zeta_1)&\psi_{W,1}(\tilde t;\zeta_1)&\ldots&\psi_{W,n-1}(\tilde t;\zeta_1)\\
					&&&\\
					\psi_{W,0}(\tilde t;\zeta_2)&\psi_{W,1}(\tilde t;\zeta_2)&\ldots&\psi_{W,n-1}(\tilde t;\zeta_2)\\
					&&&\\
					\ldots &\ldots & \ldots & \ldots\\
					&&&\\
					\psi_{W,0}(\tilde t;\zeta_n)&\psi_{W,1}(\tilde t;\zeta_n)&\ldots&\psi_{W,n-1}(\tilde t;\zeta_n)
			\end{pmatrix}$$
			where $\zeta_i$ is the $i^{th}$ root of $z$.
\end{defin}

\begin{prop}
	The matrix $\Psi_W(\tilde t;z)$ admits asympotic expansion
	$$\Psi_W(\tilde t;\Lambda)=\exp(\xi(\tilde t;\Lambda))(I+O(z^{-1}))$$
	Moreover under the isomorphism $\Xi^{-1}:H\rightarrow H^{(n)}$ we can write 
	$W\in \mathrm{Gr}^{(n)}$ as
	\begin{equation}\label{psi}
		W=\Psi_W(0,z)H^{(n)}_+
	\end{equation}
\end{prop}
\prf
	One has to note that the $i^{th}$ column of matrix $\Psi_W(\tilde t;z)$ is nothing but $\Xi^{-1}(\psi_{W,i}(\tilde t,z))$ so that asymptotic expansion follows easily.
	Equation (\ref{psi}) corresponds to the fact that $\lbrace z^{ns}\psi_{W,i}(0,z):s\in\bZ\rbrace$ is a basis for $W$.
\endprf

Observe that, since we also have
$$W=\cal W(z) H^{(n)}_+$$
we obtain 
$$\Psi_W(0,z)=\cal W(z)(I+O(z^{-1}).$$
From this equation and from lemma \ref{ToeplitzHaenkel} it follows that for every $N>0$ we have
$$T_N\Big(\cal W(\tilde t;z)\big(I+O(z^{-1})\big)\Big)=T_N(\cal W(\tilde t;z))T_N(I+O(z^{-1})).$$
Now since for every $N$
$$\det(T_N(I+O(z^{-1})))=1$$
we will assume, without loss of generality, that
$$\Psi_W(0,z)=\cal W(z)$$
since this is true modulo an irrelevant term that does not affect values of determinants we want to compute.
We now want to define a matrix $\Phi_{W}(\tilde t;z)$ analytic in $z$ near $0$ and with similar properties as $\Psi_W(\tilde t;z).$
\begin{defin}
	Let $\phi_W(\tilde t;z)$ be the unique solution of 
	$$L_W(\tilde t)\phi_W(\tilde t;z)=z^n\phi_W(\tilde t;z)$$
	analytic in $z=0$ and such that
	$$(D^i\phi)(0,z)=z^i:i=0,\ldots,n-1$$
	We define 
	$$\Phi_{W}(\tilde t;z):=\begin{pmatrix}
					1 & \zeta_1 & \ldots & \zeta_1^{n-1}\\
					&&&\\
					1 & \zeta_2 & \ldots & \zeta_2^{n-1}\\
					&&&\\
					\ldots &\ldots & \ldots & \ldots\\
					&&&\\
					1 & \zeta_n & \ldots & \zeta_n^{n-1}
			\end{pmatrix}^{-1}
			\begin{pmatrix}
					\phi_{W,0}(\tilde t;\zeta_1)&\phi_{W,1}(\tilde t;\zeta_1)&\ldots&\phi_{W,n-1}(\tilde t;\zeta_1)\\
					&&&\\
					\phi_{W,0}(\tilde t;\zeta_2)&\phi_{W,1}(\tilde t;\zeta_2)&\ldots&\phi_{W,n-1}(\tilde t;\zeta_2)\\
					&&&\\
					\ldots &\ldots & \ldots & \ldots\\
					&&&\\
					\phi_{W,0}(\tilde t;\zeta_n)&\phi_{W,1}(\tilde t;\zeta_n)&\ldots&\phi_{W,n-1}(\tilde t;\zeta_n)
			\end{pmatrix}$$
		where as before $\zeta_i$ is the $i^{th}$ root of $z$ and
		$$\phi_{W,i}(\tilde t):=D^i(\phi_W(\tilde t;z)):i=0,\ldots,n-1.$$
\end{defin}

\begin{rem}
	$\Phi_{W}(\tilde t;z)$ admits regular expansion in $z=0$ and Cauchy initial values we imposed on $\phi_W$ imply 
	$$\Phi_W(0;z)=I.$$
\end{rem}
\begin{prop}
	$\Psi_W(\tilde t;z)\Phi_W^{-1}(\tilde t;z)$ does not depend on $t_i$ for any $i$.
\end{prop}
\prf
	It is well known that equations
	$$\frac{\partial}{\partial t_i}f=(L_W^{\frac{i}{n}})_+f$$
	satisfied by $\phi_W$ and $\psi_W$ can be translated into matrix equations
	$$\frac{\partial}{\partial t_i}F=F M$$
	satisfied by $\Psi_W(\tilde t;z)$ and $\Phi_W(\tilde t;z)$
	(one can write explicitely $M$ in terms of coefficients of $(L_W^\frac{i}{n})_+$).
	Hence we have
	\begin{gather*}
		\frac{\partial}{\partial t_i}(\Psi_W(\tilde t;z)\Phi_W^{-1}(\tilde t;z))=
		\Psi_W(\tilde t;z)M\Phi_W^{-1}(\tilde t;z)-\\
		\Psi_W(\tilde t;z)\Phi_W^{-1}(\tilde t;z)\Phi_W(\tilde t;z)M\Phi_W^{-1}(\tilde t;z)=0
	\end{gather*}
\endprf

\begin{thm}\label{factorizationtheorem}
	Given a Gelfand-Dickey symbol 
	$$\cal W(\tilde t;z)=\exp\Big(\xi(\tilde t,\Lambda)\Big)\cal W(z)$$
	one can factorize it as
	$$\cal W(\tilde t;z)=\Bigg[\exp\Big(\xi(\tilde t,\Lambda)\Big)\Psi_W(-\tilde t,z)\Bigg]\Phi_W^{-1}(-\tilde t;z)$$
	where the term inside the square bracket is analytic around $z=\infty$ and the other is analytic around $z=0$.
	For assigned values of $\tilde t$ for which 
	$$\tau_{\cal W}(\tilde t)\neq 0$$
	this is the unique solution of the factorization problem (\ref{Plemelj RH}) normalized at infinity to the identity.
\end{thm}
\prf
	Using the previous proposition we have
	\begin{align*}
		\cal W(\tilde t;z)=\exp{\xi(\tilde t,\Lambda)}\cal W(z)=\exp({\xi(\tilde t,\Lambda)})\Psi_W(0,z)=\\
		\exp({\xi(\tilde t,\Lambda)})\Psi(-\tilde t;z)\Phi_W^{-1}(-\tilde t;z)\Phi_W(0;z)=\exp({\xi(\tilde t,\Lambda)})
		\Psi(-\tilde t;z)\Phi_W^{-1}(-\tilde t;z)
	\end{align*}
	Unicity of  the factorization follows from section 4.
\endprf
\begin{corol}
	For every $N>0$
	$$\tau_{W,N}(\tilde t)=\tau_{W}(\tilde t)\det(I-K_{\cal W(\tilde t;z),N})$$
	with
	\begin{equation*}
 	(K_{\cal W(\tilde t;z),N})_{ij}=\begin{cases}
 												0\quad \text{if}\:\min\lbrace i,j\rbrace< N\\
 												\\
 												\sum_{k=1}^\infty (\Psi_W(-\tilde t;z))^{(i+k)}(\Psi_W(-\tilde t;z)^{-1})^{(-j-k)}\quad\text{otherwise.}
 											\end{cases}
 	\end{equation*}
\end{corol}
\prf
	It is enough to apply Borodin-Okounkov formula using factorization obtained above.
\endprf
\begin{corol}
	For every $N>0$
	\begin{equation*}
		\dfrac{\tau_{W,N}(\tilde t)}{\tau_{W,N+1}(\tilde t)}=
		\det\Big(T\big(\Psi_W(-\tilde t;z)\big)T\big(\Psi_W(-\tilde t;z)^{-1}\big)\Big)_{N,N}
	\end{equation*}
	(observe that the right hand side of this equation is an ordinary $n\times n$ determinant, not a Fredholm determinant).
\end{corol}
\begin{prf}
	$$\dfrac{\tau_{W,N}(\tilde t)}{\tau_{W,N+1}(\tilde t)}=
	\dfrac{\det(I-K_{\cal W(\tilde t;z),N})}{\det(I-K_{\cal W(\tilde t;z),N+1})}.$$
	On the other hand the operator $(I-K_{\cal W(\tilde t;z),N+1})^{-1}(I-K_{\cal W(\tilde t;z),N})$ can be written as a block matrix obtained taking the identity matrix and replacing the $N^{th}$ block column by the $N^{th}$ block column of the matrix with $(i,j)$-entry equal to 
	$$\sum_{k=1}^\infty (\Psi_W(-\tilde t;z))^{(i+k)}(\Psi_W(-\tilde t;z)^{-1})^{(-j-k)}.$$
	Hence proof is obtained applying lemma \ref{ToeplitzHaenkel}
\end{prf}

\section{Rank one stationary reductions\\ and corresponding Gelfand-Dickey symbols}

We want to describe, more explicitely, GD symbols corresponding to solutions of Gelfand-Dickey hierarchies obtained by rank one stationary reductions.
In order to emphasize that we are dealing with rank-one generic case instead of the standard expression \emph{Krichever locus} we will speak about \emph{Burchnall-Chaundy locus}. 

\begin{defin}
	Given a point $W\in Gr^{(n)}$ we say that $W$ stays in Burchnall-Chaundy locus iff the Lax operator $L_W$ of the corresponding solution satisfies
	$$[L_W,M_W]=0$$
	for some differential operator $M_W$ of order $m$ coprime with $n$. Without loss of generality we also assume $m>n$.
\end{defin}
The name we use is due to the fact that, already in $1923$, Burchnall and Chaundy were the first to study algebras of commuting differential operators in \cite{BC} where they stated this important proposition we will use in the sequel.
\begin{prop}[\cite{BC}]
	Given a pair of commuting differential operator $L,M$ with relatively prime orders it exists an irreducible polynomial $F(x,y)$ such that
	$$F(x,y)=x^m+...\pm y^n$$
	and $F(L,M)=0.$
\end{prop}

This proposition in particular allows us to associate to every Burchnall-Chaundy solution a spectral curve defined by polynomial relation existing between the pair of commuting differential operators.
From the Grassmannian point of view one can define an action $\cal A$ of pseudodifferential operators in variable $t_1$ on $H$ by
\begin{align*}
	\cal A:\mathrm{\Psi DO}\times H&\longrightarrow H\\
	\Big((t_1)^m\frac{\partial^n}{\partial t_1^n},\varphi(z)\Big)&\longmapsto\Big(\frac{\partial^n}{\partial z^n}\Big)\Big(z^n\Big)\varphi(z) 
\end{align*}
 and, using this action, prove the following propostition
 
\begin{prop}[\cite{M}]
 	Given a point $W$ in the Burchnall-Chaundy locus one has
 	\begin{gather}
 		z^nW\subseteq W\\
 		b(z)\subseteq W
 	\end{gather}
 	where $L_W$ and $M_W$ are of order $n$ and $m$ respectively and $b(z)$ is a series in $z$ whose leading term is $z^m$. Conversely, if $W$ satisfies above properties, it stays in the Burchnall-Chaundy locus.
\end{prop}
\begin{prf}
	We just sketch the proof and make reference to Mulase's article \cite{M}.
	Suppose we are given $L_W$ and $M_W$; under conjugation with the dressing $S_W(\tilde t)$ we have
	$$S^{-1}_W(\tilde t)L_W(\tilde t)S_W(\tilde t)=\frac{\partial^n}{\partial t_1^n}$$
	Under the action $\cal A$ this gives invariance of $W$ with respect to $z^n$ while invariance with respect to $b(z)$ is obtained acting with $$S^{-1}_W(\tilde t)M_W(\tilde t)S_W(\tilde t)$$\\
	Viceversa given $W$ we reconstruct the dressing $S_W(\tilde t)$; using it we define $L_W(\tilde t)$ and $M_W(\tilde t)$ conjugating pseudodifferential operators corresponding to $z^n$ and $b(z)$.
	In particular observe that also $z^n$ and $b(z)$ will satisfy the same polynomial relation as $L_W(\tilde t)$ and $M_W(\tilde t)$.
\end{prf}

\begin{rem}
	Without loss of generality we can assume
	\begin{equation}\label{traceless}
		\frac{1}{2\pi i}\oint\frac{b(z)}{z^{ns+1}}dz=0\quad\forall s\in\bZ.
	\end{equation}
\end{rem}

Now suppose we are given an element $W=\cal W(z)H^{(n)}_+\in Gr^{(n)}$ in the Burchnall-Chaundy locus. Using the explicit isomorphism $\Xi$ we can construct a matrix $B(z):=b(\Lambda)$ such that
\begin{equation}\label{B}
	B(z)W\subseteq W.
\end{equation}

\begin{prop}\label{C}
	$$C(z):=\cal W^{-1}(z)B(z)\cal W(z)$$
	has the following properties:
	\begin{itemize}
		\item $C(z)$ is polynomial in $z$.
		\item $\mathrm{trace}(C(z))=0$
		\item $m=\max_i (j-i+n\:\mathrm{deg}\:C_{ij}(z))\quad\forall j=1\ldots n$
		\item The characteristic polynomial $p_{C(z)}(\lambda)$ of $C(z)$ defines the spectral curve of the solution.  
	\end{itemize}
\end{prop}
\begin{prf}
	Equation (\ref{B}) can be equivalently written as
	$$\cal W^{-1}(z)B(z)\cal W(z)H^{(n)}_+\subseteq H^{(n)}_+$$
	and this means precisely that $C(z)$ can't have terms in $z^{-k}$ for any $k>0$.\\
	The other properties are satisfied if and only if they are equally satisfied by $B(z)$ so that we will prove them for $B(z)$ instead of $C(z)$.
	$B(z)$ is traceless thanks to equation (\ref{traceless}) and thanks to the fact that 
	$$\mathrm{trace}(\Lambda^k)=0\quad\forall k\neq sn$$
	The third properties is satisfied as $B(z)=b(\Lambda)$ represents in $H$ multiplication by a series whose leading term is equal to $m$.
	For the last property we observe that if $F(x,y)$ is the polynomial defining the spectral curve, i.e. $F(L_W,M_W)=0$, then we will have 
	$$F(\mathrm{diag}(z,z,\ldots,z),B(z))=0$$
	as well; on the other hand thanks to Cayley-Hamilton theorem we have
	$$p_{B(z)}(B(z))=0.$$
	Since $F$ is irreducible and $p_{B(z)}(\lambda)$ has the same form
	$$p_{B(z)}(\lambda)=\lambda^n+\ldots\pm z^m$$
	we conclude that they are equal.   
\end{prf}

Observe that since $\cal W(z)$ is defined modulo multiplication on the left by invertible triangular matrices also $C(z)$ is defined modulo conjugation by elements of the group $\Delta$ of upper trianguar invertible matrices.
It was a remarkable observation of Schwarz \cite{Sc} that actually Burchnall-Chaundy locus can be described by means of matrices with properties as in proposition \ref{C} modulo the action of $\Delta$.
Here we adapt the results of \cite{Sc} to our situation.
Namely we explain how, given $C(z)$, one can recover $\cal W(z)$ and the corresponding spectral curve.

\begin{prop}\label{reconstruction}
	Given a matrix $C(z)$ such that:
	\begin{itemize}
		\item $C(z)$ is polynomial in $z$.
		\item $\mathrm{trace}(C(z))=0$
		\item $m=\max_i (j-i+n\:\mathrm{deg}\:C_{ij}(z))\quad\forall j=1\ldots n$
	\end{itemize}
	it exists a unique $W=\cal W(z)H^{(n)}_+$ in Burchnall-Chaundy locus such that its spectral curve is defined by $p_{C(z)}(\lambda)$.
\end{prop}

In order to prove this proposition we need two lemmas.

\begin{lemma}\label{C1}
	Given a polynomial matrix $C(z)$ such that
	$$m=\max_i (j-i+n\:\mathrm{deg}\:C_{ij}(z))\quad\forall j=1\ldots n$$
	(with $m$ and $n$ coprime) coefficients of characteristic polynomial
	$$p_{C(z)}(\lambda):=\lambda^n+c_1(z)\lambda^{n-1}+\ldots+c_n(z)$$ 
	satisfy
	\begin{gather*}
		n\:\mathrm{deg}\:c_s\leq ms\quad\forall s=1,\ldots,n-1\\
		\mathrm{deg}\:c_n=m
	\end{gather*}
\end{lemma}
\prf
	From 
	$$n\:\mathrm{deg}\:C_{i,j}\leq m-j+i$$ 
	and definition of determinant follows immediately that
	$$n\:\mathrm{deg}\:c_s\leq ms\quad\forall s=1,\ldots,n.$$
	Strict inequality for $s<n$ follows from the fact that $m$ and $n$ are coprime.
	For the equality
	$$\mathrm{deg}\:c_n=\mathrm{deg}\Big(\det(C(z))\Big)=m$$ 
	we observe that in every line there is a unique element $C_{ij}(z)$ such that $m=j-i+n\:\mathrm{deg}\:C_{ij}(z)$;
	taking this unique element for every line and multiplying them we will obtain the leading term of determinant which will be of order $m$.
\endprf
\begin{lemma}\label{C2}
	The equation
	\begin{equation}\label{characteristic}
			\lambda^n+c_1(z)\lambda^{n-1}+\ldots+c_n(z)=0
	\end{equation}
	with 
	\begin{gather*}
		n\:\mathrm{deg}\:c_s\leq ms\quad\forall s=1,\ldots,n-1\\
		\mathrm{deg}\:c_n=m
	\end{gather*}
	and $n,m$ coprime has $n$ distinct solutions $\lbrace \lambda_i=b(\zeta_i),\quad i=1\ldots n\rbrace$
	with 
	$$b(\zeta)=\zeta^{m}(1+O(\zeta^{-1}))$$
	(as usual $\zeta_i$ is the $i^{th}$ root of $z$).
\end{lemma}
\prf
	Imposing $\lambda_i=\zeta_i^{m}$ we have a solution of the equation
	$$(\zeta_i^m)^n+c_1(\zeta_i^n)(\zeta_i^m)^{n-1}+\ldots+c_n(\zeta_i^n)=0$$
	 at the leading order $mn$.
	Then imposing $\lambda_i=\zeta_i^{m}(1+l_1\zeta_i^{-1})$ and plugging it into the equation (\ref{characteristic}) one obtains
	$$(\zeta_i^m+l_1\zeta_i^{m-1})^n+c_1(\zeta_i^n)(\zeta_i^m+l_1\zeta_i^{m-1})^{n-1}+\ldots+c_n(\zeta_i^n)=O(\zeta_i^{mn})$$ 
	$l_1$ can be found so that terms of order $nm-1$ in the equation vanish; going on solving the equation term by term we obtain
	$$\lambda_i=\zeta_i^m\Big(1+\sum_{j<0}l_j\zeta_i^{-j}\Big)$$
	Clearly coefficients $l_j$ do not depend on the choice of the root $\zeta_i$ so that it exists $b(\lambda)$ with stated properties.
\endprf\\
Now we can prove proposition \ref{reconstruction}.\\
\prf
	We start computing the characteristic polynomial $p_{C(z)}(\lambda)$; thanks to lemmas \ref{C1} and \ref{C2} we find $n$ distinct roots $b(\zeta_1),\ldots,b(\zeta_n)$ with properties stated above.\\
	The aim is to find $\cal W(z)$ such that
	$$\cal W(z)b(\Lambda)\cal W^{-1}(z)=C(z)$$
	Since we have $n$ distinct solutions $\lbrace b(\zeta_i),\quad i=1,\ldots,n\rbrace$ of the equation
	$$p_{C(z)}(\lambda)=0$$ 
	it exists a matrix $\Upsilon(\zeta_1,\ldots,\zeta_n)$ such that	
	$$\Upsilon(\zeta_i)C(z)\Upsilon^{-1}(\zeta_i)=
	\begin{pmatrix}
		b(\zeta_1)& 0 & \ldots & 0\\
		0         & b(\zeta_2) & \ldots & 0\\
		0         & \ldots      &\ddots & 0\\
		0         & \ldots      &\ldots  &b(\zeta_n)
	\end{pmatrix}
	$$
	On the other hand it's easy to observe that the matrix $\Lambda$ can be diagonalized as
	$$\Lambda=
	\begin{pmatrix}
					1 & \zeta_1 & \ldots & \zeta_1^{n-1}\\
					1 & \zeta_2 & \ldots & \zeta_2^{n-1}\\
					\ldots &\ldots & \ldots & \ldots\\
					1 & \zeta_n & \ldots & \zeta_n^{n-1}
			\end{pmatrix}^{-1}
			\begin{pmatrix}
		\zeta_1& 0 & \ldots & 0\\
		0         & \zeta_2 & \ldots & 0\\
		0         & \ldots      &\ddots & 0\\
		0         & \ldots      &\ldots  &\zeta_n
	\end{pmatrix}
		\begin{pmatrix}
					1 & \zeta_1 & \ldots & \zeta_1^{n-1}\\
					1 & \zeta_2 & \ldots & \zeta_2^{n-1}\\
					\ldots &\ldots & \ldots & \ldots\\
					1 & \zeta_n & \ldots & \zeta_n^{n-1}
			\end{pmatrix}$$

and this means that multiplication by $b(z)$ can be written in $H^{(n)}_+$ as multiplication by
$$\begin{pmatrix}
					1 & \zeta_1 & \ldots & \zeta_1^{n-1}\\
					1 & \zeta_2 & \ldots & \zeta_2^{n-1}\\
					\ldots &\ldots & \ldots & \ldots\\
					1 & \zeta_n & \ldots & \zeta_n^{n-1}
			\end{pmatrix}^{-1}
			\begin{pmatrix}
		b(\zeta_1)& 0 & \ldots & 0\\
		0         & b(\zeta_2) & \ldots & 0\\
		0         & \ldots      &\ddots & 0\\
		0         & \ldots      &\ldots  &b(\zeta_n)
	\end{pmatrix}
		\begin{pmatrix}
					1 & \zeta_1 & \ldots & \zeta_1^{n-1}\\
					1 & \zeta_2 & \ldots & \zeta_2^{n-1}\\
					\ldots &\ldots & \ldots & \ldots\\
					1 & \zeta_n & \ldots & \zeta_n^{n-1}
			\end{pmatrix}$$
Hence we have
$$\cal W(z)=\Upsilon^{-1}(\zeta_i)\begin{pmatrix}
					1 & \zeta_1 & \ldots & \zeta_1^{n-1}\\
					1 & \zeta_2 & \ldots & \zeta_2^{n-1}\\
					\ldots &\ldots & \ldots & \ldots\\
					1 & \zeta_n & \ldots & \zeta_n^{n-1}
			\end{pmatrix}$$ 			
Note that $\cal W(z)$ is defined modulo the action of $\Delta$ so that, by construction, $C(z)$ corresponds to a unique $W\in\mathrm{Gr}^{(n)}$ such that $$W=\cal W(z)H^{(n)}_+.$$
\endprf
\begin{rem}
	As it was pointed out by Schwarz \cite{Sc}, matrices $C(z)$ with properties stated above can be used to describe points in the Grassmannian describing string solutions of Gelfand-Dickey hierarchies, i.e. solutions associated to reduction of type
	$$[L,M]=1$$
This class of solutions has not been treated in this article since they do not live in Segal-Wilson Grassmannian but just on Sato's Grassmannian constructed on the space of formal series; this means that we cannot use any more Szeg\"o-Widom theorem as  the analytical requirements are not satisfied.
Nevertheless some results obtained in section 3 still hold since the property of stability for $\lbrace\tau_{W,N}(t)\rbrace$ does not depend on analytical properties of the symbol $\cal W(z)$.
Hence one can try to apply the approach used in this article to the study of these (much less studied) string solutions; perhaps results obtained by Okounkov and Borodin in \cite{BO} for formal series and a generalization to block case can play in the setting of formal theory the same role played by Szeg\"o-Widom theorem in this paper.
\end{rem}

\begin{ex}[Symmetric $n$-coverings]\label{ex1}
	Take a symmetric $n$-covering $\cal C$ of $\bP^1$ given by equation
	\begin{equation}\label{n-covering}
		\lambda^n=P(z)=\prod_{j=1}^{nk+1}(z-a_j)
	\end{equation}
	For this particular type of curves, choosing in a appropriate way the divisor on the curve, we can write explicitely $\cal W(z),B(z)$ and $C(z)$.
	We start to observe that for any $W$ corresponding to this spectral curve we have $b(z)W\subseteq W$ with
	$$b(z)=P(z^n)^{\frac{1}{n}}$$
	Then it's easy to prove that the corresponding $B(z)=b(\Lambda)$ can be written as
	$$B(z)=\begin{pmatrix}
						0\                                    & 0\ \  & \ldots\  & 0 & z^{\frac{n-1}{n}}P(z)^{\frac{1}{n}}\\
						z^{-\frac{1}{n}}P(z)^{\frac{1}{n}}\  & 0\ \ &\ldots\ &0&0\\
						0  \                             &\ddots\ \ & \ddots\ &\vdots&\vdots\\
						\vdots \                         &\ddots\ \ &\ddots\ &\ddots&\vdots\\
						 & & & & \\
						0   \                            &\ldots\ \ &0\ &z^{-\frac{1}{n}}P(z)^{\frac{1}{n}}&0
					\end{pmatrix}$$	

	Now we define $n$ functions
	$$w_i(z):=\Big(\frac{P(z)}{z}\Big)^{\frac{i-1}{n}}\frac{1}{\prod_{j=1}^{(i-1)k}(z-a_j)},\quad i=1,\ldots,n.$$
	We take
	$$\cal W:=\mathrm{diag}(w_1(z),...,w_n(z))$$
	It is easy to verify that the matrix 
	\begin{gather*}C(z)=\cal W^{-1}(z)B(z)\cal W(z)=\\ \\ \begin{pmatrix}
						0\                                    & 0\ \  & \ldots\  & 0 & z^{\frac{n-1}{n}}P(z)^{\frac{1}{n}}\frac{ w_n(z)}{w_1(z)}\\
						z^{-\frac{1}{n}}P(z)^{\frac{1}{n}}\frac{w_1(z)}{w_2(z)}\  & 0\ \ &\ldots\ &0&0\\
						0  \                             &\ddots\ \ & \ddots\ &\vdots&\vdots\\
						\vdots \                         &\ddots\ \ &\ddots\ &\ddots&\vdots\\
						 & & & & \\
						0   \                            &\ldots\ \ &0\ &z^{-\frac{1}{n}}P(z)^{\frac{1}{n}}\frac{w_{n-1}(z)}{w_n(z)}&0
					\end{pmatrix}\end{gather*}
is polynomial in $z$.
It is worth noticing that this example already gives all possible double coverings; hence for any (possibly singular) hyperelliptc surface we found (assigning a particular divisor) the GD symbol of the corresponding algebro-geometric rank one solution of $\mathrm{KdV}$.
\end{ex}
\begin{ex}[Rational solutions]
	As pointed out by Segal and Wilson \cite{SW}, subspace of Burchnall-Chaundy locus corresponding to rational curves are given by $W=\cal W(z)H^{(n)}_+$ with $\cal W(z)$ rational in $z$.
	In particular the corresponding Gelfand-Dickey symbol will satisfy hypothesis given in proposition \ref{truncation} so that we recover the following (known) result.
	\begin{prop}
		Every rational solution of Gelfand-Dickey hierarchies can be written as a finite-size determinant.
	\end{prop}
For instance, for $n=2$, taking
$$\cal W(z)=\begin{pmatrix}
			1-d^2z^{-1} & 0\\
			0           & 1-c^2z^{-1}
		\end{pmatrix}H^{2}_+$$ 
the inverse of Gelfand-Dickey symbol is equal to 
\begin{gather*}
	\cal W^{-1}(\tilde t;z)=\\
	\begin{pmatrix}
			\cosh(z^{\frac{1}{2}}\Big(\sum_{i\geq 0}t_{2i+1}z^{2i})\Big)& 
			-z^{\frac{1}{2}}\sinh(z^{\frac{1}{2}}\Big(\sum_{i\geq 0}t_{2i+1}z^{2i})\Big)\\
			& \\
			-z^{-\frac{1}{2}}\sinh(z^{\frac{1}{2}}\Big(\sum_{i\geq 0}t_{2i+1}z^{2i})\Big)&
			\cosh(z^{\frac{1}{2}}\Big(\sum_{i\geq 0}t_{2i+1}z^{2i})\Big)
	\end{pmatrix}
	\begin{pmatrix}
			\dfrac{z}{z-d^2} & 0\\
			&\\
			0               & \dfrac{z}{z-c^2}
	\end{pmatrix}
\end{gather*}
Simply taking the residue one obtains that the corresponding $\tau$ function will be equal to
$$\tau_W(t_{1},t_3,...)=\det\begin{pmatrix}
			\cosh\Big(\sum_{i\geq 0}t_{2i+1}d^{2i+1}\Big)& 
			-d\sinh\Big(\sum_{i\geq 0}t_{2i+1}d^{2i+1}\Big)\\
			& \\
			-c^{-1}\sinh\Big(\sum_{i\geq 0}t_{2i+1}c^{2i+1}\Big)&
			\cosh\Big(\sum_{i\geq 0}t_{2i+1}c^{2i+1}\Big) 
			\end{pmatrix}$$
and recover $2$-solitons solution for $\mathrm{KdV}$.
\end{ex}
We want to point out that, for algebro geometric solutions treated in this section, the problem of factorization for Gelfand Dickey symbol can be easily translated into a Riemann-Hilbert problem on some cuts on the plane with constant jumps.
For simplicity we reduce to the case $n=2$; the procedure used here is equivalent to the one used by Its, Jin and Korepin in \cite{IJK} and generalized by Its, Mezzadri and Mo in \cite {IMM}.
Suppose we want to solve the factorization problem 
$$\cal W(\tilde t;z):=\exp\Big(\xi(\tilde t,\Lambda)\Big)\cal W(z)=T_-(\tilde t;z)T_+(\tilde t; z)$$ 
for our GD symbol with $\cal W(z)=\mathrm{diag}(w_1(z),w_2(z))$ as in example \ref{ex1}; since it will appear many times we denote $A$ the matrix
$$A:=\begin{pmatrix}
					1 & \sqrt{z}\\
					1 & -\sqrt{z}
			\end{pmatrix}$$
Also we impose 
$$P(z):=\prod_{j=1}^{2g+1}(z-a_j)$$
with all $a_j$ having modulo less then $1$ and 
$$\|a_1\|<\|a_2\|<\ldots<\|a_{2g+1}\|$$
We denote $l_1,...l_{g+1}$ the oriented intervals $(a_1,a_2),(a_3,a_4),...(a_{2g+1},\infty)$.
Instead of looking for $T_-(\tilde t;z)$ and $T_+(\tilde t;z)$ we define a new matrix $S(\tilde t;z)$ imposing
\begin{equation*}
	\begin{cases}
		$$S(\tilde t;z):=A\exp\big(-\xi(\tilde t,\Lambda)\big)T_-(\tilde t;z) \quad z\geq 1$$\\
		\\
		$$S(\tilde t;z):=A\cal W(z)T_+^{-1}(\tilde t;z)\quad z\leq 1$$
	\end{cases}
\end{equation*}
\begin{prop}\label{RHconstant}
	$S(\tilde t;z)$ has the following properties:
	\begin{itemize}
		\item It has no jumps on $S^1$
		\item It has jumps on intervals $l_j$; precisely calling $S_L(\tilde t;z)$ and $S_R(\tilde t;z)$ the values of $S(\tilde t;z)$ approaching from the left and approching from the right the interval we have
		$$S_L(\tilde t;z):=\begin{pmatrix}
													0 & 1\\
													1 & 0
												\end{pmatrix}S_R(\tilde t;z)$$
		\item It is invertible in any points but $a_j$; there it has singular behaviour of type
		$$S(\tilde t;z)\sim \begin{pmatrix}
													1 & 1\\
													1 & -1
												\end{pmatrix} (z-a_j)^{\begin{pmatrix}
													1 & 0\\
													0 & \pm\frac{1}{2}
												\end{pmatrix}}S_j(\tilde t;z)$$
												with $S_j(\tilde t;z)$ invertible in $a_j$; minus is for $a_1,\ldots,a_g$, plus for the others.
			\item At infinity it behaves as
						$$S(\tilde t;z)\sim \begin{pmatrix}
													\exp{\Big(-\sqrt{z}(t_1z+t_3z+\ldots)\Big)} & \sqrt{z}\exp{\Big(-\sqrt{z}(t_1z+t_3z+\ldots)\Big)}\\
													\exp{\Big(\sqrt{z}(t_1z+t_3z+\ldots)\Big)} & -\sqrt{z}\exp{\Big(\sqrt{z}(t_1z+t_3z+\ldots)\Big)}
												\end{pmatrix}$$
	\end{itemize}
\end{prop}

\prf
Let's call $S_+(\tilde t;z)$ and $S_-(\tilde t;z)$ the limiting values of $S(\tilde t;z)$ approaching the unit circle from inside and outside; we have
\begin{gather*}
			S_+(\tilde t;z)S_-^{-1}(\tilde t;z)=A\exp(-\xi(\tilde t,z))T_-(\tilde t;z)T_+(\tilde t;z)\cal W^{-1}(z)A^{-1}=\\
			A\exp(-\xi(\tilde t,z))\exp(\xi(\tilde t,z))\cal W(z)\cal W^{-1}(z)A^{-1}=I
\end{gather*}
and this proves we haven't any jumps on $S^1$.\\
Writing explicitely $S(\tilde t;z)$ as 
\begin{equation*}
	\begin{cases}
		$$S(\tilde t;z)=
		\begin{pmatrix}
			\exp{\Big(-\sqrt{z}(t_1z+t_3z+\ldots)\Big)} & \sqrt{z}\exp{\Big(-\sqrt{z}(t_1z+t_3z+\ldots)\Big)}\\
			&\\
			\exp{\Big(\sqrt{z}(t_1z+t_3z+\ldots)\Big)} & -\sqrt{z}\exp{\Big(\sqrt{z}(t_1z+t_3z+\ldots)\Big)}
		\end{pmatrix}T_-(\tilde t;z) \quad z\geq 1$$\\
		\\
		$$S(\tilde t;z)=\begin{pmatrix}
												1 & \dfrac{(P(z))^{\frac{1}{2}}}{\prod_{j=0}^g(z-a_j)}\\
												&\\
												1 & -\dfrac{(P(z))^{\frac{1}{2}}}{\prod_{j=0}^g(z-a_j)}
											\end{pmatrix}T_+^{-1}(\tilde t;z)\quad z\leq 1$$
	\end{cases}
\end{equation*}
we obtain almost immediately the other points of the proposition; the only thing we have to observe is that both $T_+(\tilde t;z)$ and $T_-(\tilde t;z)$ are invertible inside and outside the circle respectively.
This is because we have
$$\det{\cal W(\tilde t;z)}=\dfrac{(P(z))^{\frac{1}{2}}}{\prod_{j=0}^g(z-a_j)}=\det(T_+(\tilde t;z))\det(T_-(\tilde t;z))$$
This condition combined with 
$$\lim _{z\rightarrow\infty}\det(T_-(\tilde t;z))=1$$
gives 
\begin{gather*}
	\det(T_+(\tilde t;z))=1\\
	\det(T_-(\tilde t;z))=\dfrac{(P(z))^{\frac{1}{2}}}{\prod_{j=0}^g(z-a_j)}.
\end{gather*}
\endprf\\
The Riemann-Hilbert problem given by propostion \ref{RHconstant} is equivalent to the one proposed in section 5.
What can be done is to write explicitely the solution $S(\tilde t;z)$ using $\theta$ functions associated to the curve; this is what has been done in \cite{IJK} and \cite{IMM}. Actually comparing previous proposition with results obtained in section 5 we immediately realize that 
$$S(\tilde t;z)=\begin{pmatrix}
					\psi_{W,0}(\tilde t;\sqrt{z})&\psi_{W,1}(\tilde t;\sqrt{z})\\
					& \\
					\psi_{W,0}(\tilde t;-\sqrt{z})&\psi_{W,1}(\tilde t;-\sqrt{z})\\
			\end{pmatrix}$$
so that all we have to do in our case is to write down Baker-Akhiezer function in terms of special functions.
We can carry on the same procedure for $n$ arbitary; the only difference will be that the jump matrices will remain constant but more complicated; in any case the solution of this Riemann-Hilbert problem with constant jumps will be
$$S(\tilde t;z)=\begin{pmatrix}
					\psi_{W,0}(\tilde t;\zeta_1)&\psi_{W,1}(\tilde t;\zeta_1)&\ldots&\psi_{W,n-1}(\tilde t;\zeta_1)\\
					&&&\\
					\psi_{W,0}(\tilde t;\zeta_2)&\psi_{W,1}(\tilde t;\zeta_2)&\ldots&\psi_{W,n-1}(\tilde t;\zeta_2)\\
					&&&\\
					\ldots &\ldots & \ldots & \ldots\\
					&&&\\
					\psi_{W,0}(\tilde t;\zeta_n)&\psi_{W,1}(\tilde t;\zeta_n)&\ldots&\psi_{W,n-1}(\tilde t;\zeta_n)
			\end{pmatrix}$$
Explicit formulas involving special functions can be used here to apply proposition \ref{factorization} to our case.
For instance taking the elliptic curve $\cal C$ given by equation
$$w^2=4 z^3-g_2z-g_3$$ 
with uniformization given by the Weierstass $\wp$ function 
$$(z,w)=(\wp(u),\wp'(u))$$ 
one can write wave function as
$$\psi(x,t,u):=\frac{\sigma(u-c-x)\sigma(c)}{\sigma(u-c)\sigma(x+c)}\exp\Big(x\zeta(u)-\unmez t\wp'(u)\Big)$$
(here $\zeta$ and $\sigma$ are Weierstrass $\zeta$ and $\sigma$ function respectively, $x$ and $t$ correspond to the first and the third time).
With some tedious computations, making the change of variables $u=u(z)$, the right hand side of equation (\ref{RH Widom}) can be obtained.
It turns out that the only relevant factorization is the one given by
$$\cal W^{-1}(x,t;u)=\Big[\cal W^{-1}(u)\Psi(-x,-t;u)\Big]\Big[\Psi^{-1}(-x,-t;u)\exp(-x\Lambda-t\Lambda^3)\Big]$$
where as before (we just wrote $z$ as a function of $u$) we have
$$\Psi(x,t;u):=\begin{pmatrix}
           1 & (\wp(u))^{\unmez}\\
		   1 & -(\wp(u))^{\unmez}
		\end{pmatrix}^{-1}\begin{pmatrix}
           \psi(x,t,u) & \partial _x\psi(x,t,u)\\
		   \psi(x,t,-u) & \partial _x\psi(x,t,-u)
		\end{pmatrix}$$
Plugging into equation (\ref{RH Widom}) we obtain
$$\frac{d}{dx}\tau(x,t)=K t+2\zeta(-c)-2\zeta(x-c)$$
(here $K$ is some constant); taking another derivative we obtain elliptic solution of $\mathrm{KdV}$ as expected.

\end{document}